\documentclass[preprint]{elsarticle}

\journal{Computer Methods in Applied Mechanics and Engineering}

\usepackage{lipsum}
\usepackage{amsfonts}
\usepackage{graphicx}
\usepackage{epstopdf}
\usepackage{algorithm}
\usepackage{algorithmic}
\usepackage{amsmath}
\ifpdf
  \DeclareGraphicsExtensions{.eps,.pdf,.png,.jpg}
\else
  \DeclareGraphicsExtensions{.eps}
\fi
\usepackage{float}
\usepackage{graphicx}
\usepackage{booktabs} 
\usepackage[table]{xcolor}
\usepackage{bm}
\usepackage{setspace}
\usepackage{enumitem}
\usepackage{graphicx}
\usepackage[mathscr]{euscript}
\usepackage{cite}

\newtheorem{assumption}{Assumption}

\renewcommand{\vec}[1] {\ensuremath{\boldsymbol{#1}}}

\usepackage{amsopn}

\def\R{{\mathbb R}}
\def\U{{\mathcal U}}
\def\Z{{\mathcal Z}}
\def\t{{\vec{\theta}}}
\def\u{{\vec{u}}}

\def\z{{\vec{z}}}
\def\zbar{{\overline{\z}}}
\def\x{{\vec{x}}}
\def\e{{\vec{e}}}

\def\y{{\vec{y}}}

\def\b{{\vec{b}}}
\def\J{{\vec{J}}}

\def\d{{\vec{\delta}}}
\def\A{{\vec{A}}}
\def\B{{\vec{B}}}

\def\M{{\vec{M}}}

\def\G{{\vec{G}}}
\def\S{{\vec{S}}}

\def\L{{\vec{L}}}
\def\C{{\vec{C}}}

\def\H{{\vec{H}}}

\def\X{{\vec{X}}}

\def\I{{\vec{I}}}
\def\g{{\vec{g}}}

\def\F{{\vec{F}}}
\DeclareMathAlphabet\mathbfcal{OMS}{cmsy}{b}{n}

\usepackage{array}
\newcolumntype{P}[1]{>{\centering\arraybackslash}p{#1}}
\newcolumntype{M}[1]{>{\centering\arraybackslash}m{#1}}

\usepackage[a4paper, total={6in, 9in}]{geometry}

\begin{document}

\begin{frontmatter}
  
\title{Hyper-differential sensitivity analysis with respect to model discrepancy: optimal solution updating}

\author[1]{Joseph Hart%
   \fnref{fn1}}
\ead{johart@sandia.gov}
\author[2]{Bart van Bloemen Waanders\fnref{fn2}\corref{cor1}}
\ead{bartv@sandia.gov}

\cortext[cor1]{Corresponding author}

\affiliation[1]{organization={Sandia National Laboratories},
  addressline={P.O. Box 5800, Albuquerque, NM 87185 },
  country={United States}}

\affiliation[2]{organization={Sandia National Laboratories},
  addressline={P.O. Box 5800, Albuquerque, NM 87185 },
 country={United States}}

\begin{abstract}
  A common goal throughout science and engineering is to solve optimization problems constrained by computational models. However, in many cases a high-fidelity numerical emulation of systems cannot be optimized due to code complexity and computational costs which prohibit the use of intrusive and many query algorithms. Rather, lower-fidelity models are constructed to enable intrusive algorithms for large-scale optimization. As a result of the discrepancy between  high and low-fidelity models, optimal solutions determined using  low-fidelity models are frequently far from true optimality. In this article we introduce a novel approach that uses post-optimality sensitivities with respect to model discrepancy  to update the optimization solution. Limited high-fidelity data is used to calibrate the model discrepancy in a Bayesian framework which in turn is propagated through post-optimality sensitivities of the low-fidelity optimization problem. Our formulation exploits structure in the post-optimality sensitivity operator to achieve computational scalability. Numerical results demonstrate how an optimal solution computed using a low-fidelity model may be significantly improved with limited evaluations of a high-fidelity model. 
\end{abstract}


\end{frontmatter}

\section{Introduction}

Optimization problems constrained by computational models are ubiquitous in science and engineering and have been the topic of extensive research \citep{Vogel_99, Archer_01,Haber_01,Vogel_02,Biegler_03,Biros_05,Laird_05,Hintermuller_05,Hazra_06,Biegler_07,Borzi_07,Hinze_09,Biegler_11,frontier_in_pdeco}. Nonetheless, optimization problems are only as useful as the models that constrain them. In many applications, high-fidelity models cannot be instrumented for optimization due to their software complexity and computational cost. This is common in production codes which are developed over many years with the goal of achieving model fidelity rather than enabling optimization or other intrusive analysis approaches. In such cases, a logical option is to construct low-fidelity models and thereby enable efficient optimization.  Assuming that the optimization problem can be solved using a low-fidelity model, we pose the question: ``How can a small number of evaluations from the high-fidelity model be used to improve the optimal solution?". This would benefit a variety of complex optimization problems for applications with multi-scale, multi-physics, nonlinear, and multi-component features (such as climate, additive manufacturing, and experimental fusion).  

In this article we introduce a novel extension of hyper-differential sensitivity analysis (HDSA) to compute the sensitivity of optimization problems with respect to model discrepancy, i.e. the difference in the high and low-fidelity models. HDSA is built on the foundation of post-optimality sensitivity analysis~\citep{shapiro_SIAM_review,fiacco} which was originally developed in the context of operations research and then extended to optimization constrained by partial differential equations (PDEs) targeting stability analysis~\citep{Griesse_part_1,Griesse_part_2,griesse2}. HDSA was specifically created to evaluate the influence of uncertainties by scaling post-optimality sensitivities to high-dimensions through a coupling of tools from PDE-constrained optimization and numerical linear algebra~\citep{HDSA,saibaba_gsvd,sunseri_hdsa,hart_2021_bayes,sunseri_2,hart_hdsa_control_oed}. 

Our extension of HDSA provides a tool to predict how the optimal solution will change given a perturbation of the model discrepancy. However, since the space of possible model discrepancies is large (the discrepancy is an operator between function spaces), it is critical to constrain the sensitivity calculation by high-fidelity data. Accordingly we formulate a Bayesian inverse problem to estimate the model discrepancy and compute sensitivities to update the optimal solution in the direction of the maximum a posterior probability point. A systematic integration of prior knowledge, high-fidelity data, and low-fidelity optimization is enabled in a pragmatic way for large-scale applications.

Calibration of model discrepancy in a Bayesian framework was proposed in the seminal work of Kennedy and O'Hagan~\citep{ohagan2001} to account for uncertainty and correlations between the discrepancy and model parameters being estimated. This was followed by many works which analyzed and extended their approach,~\citep{ohagan2002,ohagan2014,Higdon_2008,Arendt_2012,Maupin,Ling_2014,Gardner_2021}.  In the calibration context, the model discrepancy is accessed through observed data with an implicit dependence on the calibration parameters which are not known a priori.  Our work is inspired by these concepts but is significantly different in our assumption that a high-fidelity model may be queried for different values of the optimization variables as specified by the user.  We focus on optimal control and design problems, although model calibration may still be considered from a different perspective than Kennedy and O'Hagan.

Our motivation to incorporate limited high-fidelity model evaluations to improve analysis using a low-fidelity model is not new. Multi-fidelity methods are motivated by a need to enable outer loop analysis such as optimization and uncertainty quantification by combining evaluations of high and low-fidelity models to achieve accuracy commensurate to high-fidelity analysis using a smaller number of high-fidelity model evaluations~\citep{multifidelity_review_peherstorfer}.  There has been considerable interest in the context of accelerating Monte Carlo type analysis~\citep{multifidelity_gorodetsky} for general uncertainty quantification applications as well as in the specific context of optimization under uncertainty~\citep{multi_fidelity_opt_willcox,multifidelity_wind_plant_geraci}. There has also been work to combine high and low-fidelity model evaluations within each iterate of a deterministic optimization algorithm~\citep{multifidelity_quasinewton_bryson,trmm_lewis,biegler_rom_opt,willcox_multifi_opt_2012}. Analogously, we assume access to both a high and low-fidelity model.  However, unlike the previous references we seek a framework where the optimization problem does not require access to the high-fidelity model. We also assume that derivatives of the high-fidelity model are not available. Furthermore, our focus is on cases where the number of high-fidelity solves is significantly constrained and must be performed offline. We cannot provide theoretical guarantees which multi-fidelity methods ensure through recourse on the high-fidelity model. Our framework is suited for applications where multi-fidelity optimization is infeasible due to the complexity of the high-fidelity model. 

Our contributions in this article consist of:
\begin{enumerate}
\item[$\bullet$] In section~\ref{sec:post_opt_sen} we provide a formulation and derivation to compute post-optimality sensitivities with respect to model discrepancy. Furthermore, we enable computational scalability by exploiting structure in the post-optimality sensitivities to constraint the form of the model discrepancy.
\item[$\bullet$] In section~\ref{sec:method_and_derivation} we pose a Bayesian inverse problem to estimate the discrepancy (which is an operator) from high-fidelity data, determine a computationally efficient expression for the posterior discrepancy, and propagate it through the post-optimality sensitivities.
\item[$\bullet$] In sections~\ref{sec:prior_discrepancy} and~\ref{sec:numerics} we present an approach to specify the prior and noise model.
\item[$\bullet$] In section~\ref{sec:numerics} we demonstrate the proposed approach on an illustrative example and a more complex fluid flow example where the Stokes equation is used as a low-fidelity approximation of the Naiver-Stokes equation.
\end{enumerate}

\section{Optimization Formulation} \label{sec:optimization_formulation}
Consider optimization problems of the form
\begin{align}
\label{eqn:true_opt_prob}
& \min_{z \in \Z} J(S(z),z) 
\end{align}
where $z$ denotes optimization variables in the Hilbert space $\Z$, $S:\Z \to \U$ denotes the solution operator for a PDE with state variable $u$ in an infinite dimensional Hilbert space $\U$, and $J:\U \times \Z \to \R$ is the objective function. This formulation is applicable to design, control, and inverse problems. For simplicity we will refer to $z$ as the controller, but emphasize that our formulation is general. We focus on problems where $\Z$ is infinite dimensional.

In many applications, evaluating $S(z)$, i.e. solving the PDE, is computationally intensive, so in practice we solve
\begin{align}
\label{eqn:approx_opt_prob}
& \min_{z \in \Z} J(\tilde{S}(z),z) 
\end{align}
where $\tilde{S}:\Z \to \U$ is the solution operator for a lower-fidelity PDE. Common practice is to solve~\eqref{eqn:approx_opt_prob} and use limited evaluations of $S(z)$ to assess the quality of the solution. Our goal in this article is to use the high-fidelity model evaluations to improve the solution computed in~\eqref{eqn:approx_opt_prob}.

\section{Post-optimality sensitivity with respect to model discrepancy} \label{sec:post_opt_sen}
Our derivations will leverage the PDE discretization to define the model discrepancy. To this end, we briefly introduce the essential components of the discretization. Let $\{\phi_i\}_{i=1}^m$ be a basis for a finite dimensional subspace $\U_h \subset \U$ and $\{\psi_j\}_{j=1}^n$ be a basis for a finite dimensional subspace of $\Z_h \subset \Z$. These may, for instance, be finite element basis functions. Let $\u \in \R^m$ and $\z \in \R^n$ denote coordinates of $u$ and $z$, respectively, in these bases. The solution operators are discretized by $\S:\R^n \to \R^m$ and $\tilde{\S}:\R^n \to \R^m$ and the objective function is discretized by $\J:\R^m \times \R^n \to \R$.

To propagate model discrepancy through the optimization problem, let $\d(\z,\t)$ be a function of $\z$, parameterized by $\t$ (which will be defined below), and consider
\begin{align}
\label{eqn:dis_approx_opt_prob_pert_rs}
 \min_{\z \in \R^n} \hspace{1 mm} \widehat{\J}(\z,\t):=\J(\tilde{\S}(\z)+\d(\z,\t),\z) .
\end{align}
The parameterized optimization problem~\eqref{eqn:dis_approx_opt_prob_pert_rs} coincides with the low-fidelity optimization problem when $\d \equiv \vec{0}$ and the high-fidelity optimization problem when $\d=\S-\tilde{\S}$. Let $\t_{\text{nom}}=\vec{0}$ and assume that $\d$ is parameterized such that $\d(\z,\t_{\text{nom}})=\vec{0}$ $\forall \z$. 

To facilitate our subsequent analysis, we assume that $\widehat{\J}$ is twice continuously differentiable with respect to $(\z,\t)$. Let $\tilde{\z} \in \R^n$ denote a local minimum of the low-fidelity problem,  i.e. a minimizer of the parameterized problem~\eqref{eqn:dis_approx_opt_prob_pert_rs} when $\t=\t_\text{nom}$. We assume that $\tilde{\z}$ satisfies the first and second order optimality conditions, 
\begin{align}
\label{eqn:first_order_opt}
\nabla_\z \widehat{\J}(\tilde{\z},\t_{\text{nom}})=\vec{0}
\end{align}
and $\nabla_{\z,\z} \widehat{\J}(\tilde{\z},\t_{\text{nom}})$ is positive definite, respectively, where $\nabla_z$ and $\nabla_{\z,\z}$ denote the gradient and Hessian with respect to $\z$, respectively.

Applying the Implicit Function Theorem to the first order optimality condition~\eqref{eqn:first_order_opt}, we consider an operator $\F: \mathcal N(\t_{\text{nom}}) \to \mathcal N(\tilde{\z})$, defined on neighborhoods of $\t_{\text{nom}}$ and $\tilde{\z}$, such that $\F(\t)$ is a stationary point of~\eqref{eqn:dis_approx_opt_prob_pert_rs}, i.e.
$$\nabla_\z \widehat{\J}(\F(\t),\t)=\vec{0} \qquad \forall \t \in \mathcal N(\t_{\text{nom}}).$$
Assuming that $\nabla_{\z,\z} \widehat{\J}(\F(\t),\t)$ is positive definite for all $\t \in \mathcal N(\t_{\text{nom}})$, $\F$ maps the model discrepancy to the optimal solution. The Jacobian of $\F$ with respect to $\t$, i.e. the sensitivity of the optimal solution to the model discrepancy, evaluated at $\t=\vec{0}$, is given by
\begin{eqnarray}
\label{eqn:dis_sen_op}
\nabla_\t \F(\vec{0}) = - \H^{-1} \B \in \R^{n \times p}
\end{eqnarray}
where $\H =\nabla_{\z,\z} \widehat{\J}(\tilde{\z},\vec{0}) \in \R^{n \times n}$ is the Hessian of $\widehat{\J}$ and $\B = \nabla_{\z,\t} \widehat{\J}(\tilde{\z},\vec{0}) \in \R^{n \times p}$ is the Jacobian of $\nabla_\z \widehat{\J}$ with respect to $\t$.  Both are evaluated at the low-fidelity solution $(\tilde{\z},\vec{0})$. We interpret the Jacobian matrix $\nabla_{\t} \F(\vec{0})$ as a Newton update of the optimal solution $\tilde{\z}$ when the model discrepancy is perturbed.

To determine an appropriate parameterization of $\d(\z,\t)$, we consider (in the function spaces) a general operator mapping from $\mathcal Z_h$ to $\mathcal U_h$. Such a operator can be written as
\begin{align*}
&z \mapsto \sum\limits_{i=1}^m f_i(z) \phi_i
\end{align*}
where $f_i:\Z_h \to \R$, $i=1,2,\dots,m$, denotes functionals, which will be parameterized by $\t$. We realize a simplification by observing that the sensitivity operator~\eqref{eqn:dis_sen_op} only depends on $(\z,\z)$ and $(\z,\t)$ derivatives of the objective $\widehat{\J}$ evaluated at $(\tilde{\z},\vec{0})$, so without loss of generality we assume that $f_i(z)$, $i=1,2,\dots,m$, are affine functions of $z$. Using the Riesz representation for bounded linear functionals we have a general form 
\begin{eqnarray*}
 z \mapsto \sum\limits_{i=1}^m (\theta_{i,0}+(z,a_i(\t))_{\mathcal Z}) \phi_i
\end{eqnarray*}
where $\theta_{i,0} \in \R$, $i=1,2,\dots,m$ and $a_i(\t) \in Z_h$, $i=1,2,\dots,m$, are the Riesz representation elements, which we parameterize with $\t$. A general expression for elements in $Z_h$ is given by writing them as a linear combination of basis functions as
\begin{eqnarray*}
a_i(\t) =  \sum\limits_{j=1}^n  \theta_{i,j} \psi_j,
\end{eqnarray*}
where $\theta_{i,j} \in \R$. This yields an operator of the form
\begin{align*}
&z \mapsto \sum\limits_{i=1}^m \left( \theta_{i,0} + \sum\limits_{j=1}^n \theta_{i,j} (z,\psi_j)_{\mathcal Z} \right)\phi_i
\end{align*}
where the vector of coefficients is defined as $\t=(\t_0^T,\t_1^T,\dots,\t_m^T)^T \in \R^p$, $p=m(n+1)$, where $\t_0=(\theta_{0,0},\theta_{1,0},\dots,\theta_{m,0})^T \in \R^m$ corresponds to the $m$ intercept terms and $\t_i = (\theta_{i,1},\theta_{i,2},\dots,\theta_{i,n})^T \in \R^n$, $i=1,2,\dots,m$, corresponds to the $m$ linear functionals.

Transforming into the coordinate spaces, the model discrepancy $\d:\R^n \times \R^p \to \R^m$ can be written in a convenient Kronecker product representation 
\begin{eqnarray}
\label{eqn:delta_kron}
\d(\z,\t) =
\left( \begin{array}{cc}
\I_m & \I_m \otimes \z^T \M_z
\end{array} \right) \t
\end{eqnarray}
where $\I_m \in \R^{m \times m}$ is the identity matrix and $\M_z \in \R^{n \times n}$ is the mass matrix whose $(i,j)$ entry is $(\psi_i,\psi_j)_Z$. The Kronecker product structure of~\eqref{eqn:delta_kron} will prove critical in the subsequent analysis to achieve computational scalability as $p$, the dimension of $\t$, is extremely large.

The post-optimality sensitivity operator~\eqref{eqn:dis_sen_op} depends on two matrices, the inverse Hessian $\H^{-1}$ and the mixed second derivative $\B$. The Hessian $\H$ does not depend on our representation of the discrepancy since it is evaluated at $\t=\vec{0}$ and $\d(\z,\vec{0}) \equiv \vec{0}$. To determine $\B = \nabla_{\z, \t} \widehat{\J}(\tilde{\z},\vec{0})$ we apply the Chain and Product rules to differentiate $\widehat{\J}(\z,\t) = \J(\tilde{\S}(\z)+\d(\z,\t),\z)$. This gives  
\begin{eqnarray}
\label{eqn:B_abstract}
\nabla_{\z,\t} \widehat{\J} = \nabla_\z\tilde{\S}^T \nabla_{\u,\u} \J \nabla_\t\d + \nabla_\z\d^T \nabla_{\u,\u} J \nabla_\t\d + \nabla_\u \J \nabla_{\z,\t} \d + \nabla_{\z,\u} \J \nabla_\t\d ,
\end{eqnarray}
where we adopt the convention that gradients are row vectors. 

To simplify the subsequent analysis we assume\footnote{This is common on a wide range of problems when the objective function admits the structure $J(u,z) = J_{mis}(u) + J_{reg}(z)$ where $J_{mis}$ is a state misfit or design criteria and $J_{reg}$ is an optimization variable regularization.} that $ \nabla_{\z,\u} \J =0$. Our derivations can be done in the more general case but such additional algebraic manipulations do not add value to the exposition. Observe that $\nabla_\z \d(\tilde{\z},\vec{0})=\vec{0}$ since $\d(\z,\vec{0}) \equiv \vec{0}$. The remaining nonzero derivatives of $\d$ in~\eqref{eqn:B_abstract} can be computed thanks to the form of the discrepancy~\eqref{eqn:delta_kron}. They are given by
\begin{eqnarray*}
\nabla_\t \d(\z,\t) =
\left( \begin{array}{cc}
 \I_m &   \I_m \otimes \z^T \M_z 
 \end{array} \right)
 \qquad \text{and} \qquad
 \nabla_\u \J \nabla_{\z,\t}\d(\z,\t) =
\left( \begin{array}{cc}
\vec{0} & \nabla_\u \J  \otimes \M_z 
\end{array} \right).
\end{eqnarray*}
It follows that
\begin{eqnarray}
\label{eqn:B}
\B = 
\nabla_\z \tilde S^T \nabla_{\u, \u} \J
\left( \begin{array}{cc}
 \I_m & \I_m \otimes \tilde{\z}^T \M_z 
 \end{array} \right)
 +
\left( \begin{array}{cc}
\vec{0} & \nabla_\u \J \otimes \M_z
 \end{array} \right)
 \in \R^{n \times p} 
\end{eqnarray}
where $\nabla_\z \tilde \S$, $\nabla_\u \J$, and $\nabla_{\u,\u} \J$ are evaluated at $(\tilde{\S}(\tilde{\z}),\tilde{\z})$. We emphasize that the Kronecker product structure of $\d$ is preserved in $\B$.
\section{Updating the optimal solution} \label{sec:method_and_derivation}
Assume that the low-fidelity optimization problem has been solved to determine $\tilde{\z}$, and that we have access to $N$ forward solves of the high-fidelity model $\S$ for different inputs $\{ \z_\ell \}_{\ell=1}^N$. Let $\y_\ell=\S(\z_\ell)-\tilde{\S}(\z_\ell)$, $\ell=1,2,\dots,N,$ denote the corresponding evaluations of the model discrepancy at these inputs. The number of evaluations $N$ will be small since evaluating $\S(\z_\ell)$ is computationally costly. Our goal is to use these limited high-fidelity evaluations to improve the optimal solution.

If $N<<n$, the number of high-fidelity solves is much less than the dimension of $\z$, then estimating the model discrepancy is ill-posed. We formulate a Bayesian inverse problem to incorporate domain expertise (prior knowledge) and uncertainty quantification into our estimation of the model discrepancy parameters $\t$. We use $\vec{\Theta}$ to denote a random vector taking values in $\R^p$ which models uncertainty in the discrepancy parameters $\t$ and $\d(\z,\vec{\Theta})$ to denote the stochastic operator which corresponds to the push forward of $\vec{\Theta}$ through the deterministic operator $\d(\z,\t)$. We likewise consider the push forward of $\vec{\Theta}$ through the post-optimality sensitivity operator, shifted by the low-fidelity optimal solution, 
\begin{align*}
\vec{U} = \tilde{\z} + \nabla_{\t} \F(\tilde{\z},\vec{0}) \vec{\Theta}.
\end{align*}
The random vector $\vec{U}$ takes values in $\R^n$. The mean of $\vec{U}$ is our updated optimal solution given the high-fidelity data. This article focuses on computing the posterior mean of $\vec{\Theta}$ and propagating it through the post-optimality sensitivities to compute the posterior mean of $\vec{U}$. A follow on article will explore posterior sampling to quantify uncertainty. 

\begin{figure}[h]
\centering
  \includegraphics[width=0.95\textwidth]{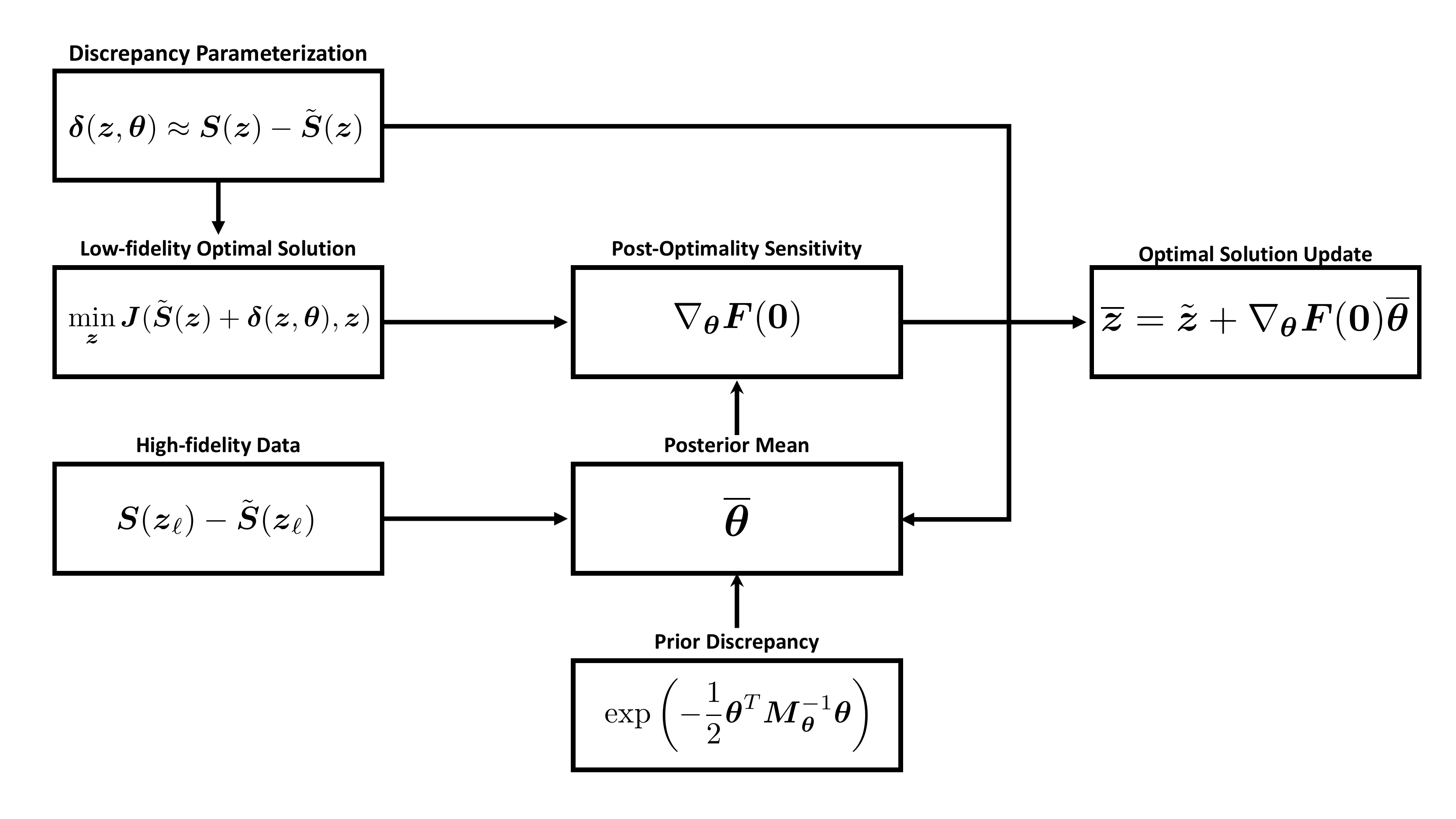}
  \caption{Diagram of the proposed analysis pipeline. We denote the discrepancy prior covariance and posterior mean as $\M_\t^{-1}$ and $\overline{\vec{\theta}}$, respectively. $\tilde{\z}$ and $\zbar$ denotes the low-fidelity and updated optimal solutions, respectively.}
  \label{fig:algorithm_diagram}
\end{figure}

Figure~\ref{fig:algorithm_diagram} depicts our proposed approach to combine Bayesian inversion and post-optimality sensitivities to update the optimal solution. In this section, we detail our contributions to enable this workflow, which include:
\begin{enumerate}
\item[1.] defining a prior distribution for $\vec{\Theta}$,
\item[2.] formulating a Bayesian inverse problem to estimate the discrepancy,
\item[3.] determining an efficient closed form expression for the posterior mean of $\vec{\Theta}$, and
\item[4.] propagating the posterior mean through the post-optimality sensitivity operator to compute the mean of $\vec{U}$.
\end{enumerate}

Due to the high dimensionality of the model discrepancy parameterization, $\t \in \R^p$, it is critical that the Kronecker product structure of the discrepancy is preserved in the posterior of $\vec{\Theta}$.  We make judicious choices in the problem formulation and manipulate the linear algebra to achieve this. The approach never requires computation in $\R^p$, and hence is efficient for large-scale optimization problems. 

\subsection{Defining the prior}
The parameters $\vec{\theta}$ correspond to the discretization of an operator mapping between infinite dimensional function spaces. Infinite dimensional Bayesian inverse problems~\citep{ghattas_infinite_dim_bayes_1,ghattas_infinite_dim_bayes_2} have seen great advances in recent years. However, most existing research focuses on estimating functions rather than operators (a mapping between functions). This subsection details our approach to define a prior on the space of operators by exploiting the Kronecker product structure of $\d$, which is crucial to ensure computational scalability.

To facilitate computation and ease specification of prior parameters, we take a Gaussian prior for $\vec{\Theta}$ with mean $\vec{0}$ and a covariance matrix to incorporate application knowledge into the discrepancy. Such priors are common in infinite dimensional Bayesian inverse problems. Recall that a mean zero Gaussian random vector $\vec{X}$ with covariance $\vec{\Pi}$ has a probability density function proportional to $\exp\left( - \frac{1}{2} \vec{x} ^T \vec{\Pi}^{-1} \vec{x}  \right)$. The highest probability regions correspond to the points $\vec{x}$ for which the weighted norm
$ \vert \vert \vec{x}  \vert \vert_{\vec{\Pi}^{-1}}=\sqrt{ \vec{x}^T \vec{\Pi}^{-1} \vec{x} }$ is small. Hence, we define a prior covariance for $\vec{\Theta}$ by determining a symmetric positive definite matrix which defines a norm on the space of operators.

Since $\t \in \R^p$ does not have a physical interpretation, we start by defining a prior in $\R^m$ (or equivalently $\U_h$). Let $\L \in \R^{m \times m}$ be a symmetric positive definite matrix which defines a weighted norm on the state space for which small norm elements will conform to discrepancy characteristics specified from domain expertise. This may embed assumptions about the smoothness of $\d$, boundary conditions, or penalize deviations from known conservation properties. A common approach is to define $\L$ as the square of an elliptic differential operator so that $\L^{-1}$ corresponds to the discretization of a well-defined covariance operator in the infinite dimensional state space~\citep{ghattas_infinite_dim_bayes_1,ghattas_infinite_dim_bayes_2,stuart_inv_prob}.

Computing the $\L$-weighted inner product of $\d(\z,\t)$ with itself and exploiting Kronecker product structure we have
\begin{align*}
( \d(\z,\t),\d(\z,\t))_\L = \t^T 
 \left( \begin{array}{cc}
\L & \L \otimes \z^T \M_z \\
\L \otimes \M_z \z & \L \otimes \M_\z \z \z^T \M_\z 
\end{array} \right)
\t .
\end{align*}
To remove dependence on $\z$, we integrate $( \d(\z,\t),\d(\z,\t))_\L$ with respect to a Gaussian measure on $\R^n$ with mean $\tilde{\z}$ and covariance matrix $\zeta^2 \M_\z^{-1}$. This measures the average size of the model discrepancy giving the greatest weight to the low-fidelity solution $\tilde{\z}$ and decreases the weight proportionally (with constant $\zeta$) with the controllers' distance (in $\Z_h$) from $\tilde{\z}$. Recalling the property of Gaussians that the expectation of $\z \z^T$ is $\zeta^2 \M_\z^{-1} + \tilde{\z} \tilde{\z}^T$, and simplifying algebra, yields the symmetric positive definite matrix 
\begin{eqnarray}
\label{eqn:M_theta}
\M_\t =  \left( \begin{array}{cc}
\L & \L \otimes \tilde{\z}^T \M_z \\
\L \otimes \M_z \tilde{\z} & \L \otimes (\zeta^2 \M_\z+ (\M_\z \tilde{\z}) (\M_\z \tilde{\z})^T)
\end{array} \right)  \in \R^{p \times p},
\end{eqnarray} 
for which
\begin{eqnarray}
\label{eqn:M_theta_ip}
\t^T \M_\t \t = \int_{\R^n} \vert \vert \d(\z,\t) \vert \vert_{\L}^2 d \mu(\z),
\end{eqnarray}
where $\mu$ is the Gaussian measure on $\R^n$ with mean $\tilde{\z}$ and covariance $\zeta^2 \M_\z^{-1}$.

Hence $\M_\t$ defines an inner product for $\t$ to measure the size of the model discrepancy $\d(\z,\t)$ according to our prior knowledge imposed by $\L$. Let $\M_\t^{-1}$ be the prior covariance for $\vec{\Theta}$, which ensures that the prior is rooted in established theory relating known physical properties to covariance matrices which can be efficiently manipulated in the state space.

As a building block for our subsequent analysis, note that we can decompose $\M_\t = \vec{C} \vec{C}^T$ where
\begin{align}
\label{eqn:C}
\vec{C} =
\left(
\begin{array}{cc}
\L^{\frac{1}{2}} & \vec{0} \\
\L^{\frac{1}{2}}  \otimes \M_\z \tilde{\z} & \L^{\frac{1}{2}}  \otimes \zeta \M_\z^{\frac{1}{2}} 
\end{array}
\right) 
\end{align}
and 
\begin{align}
\label{eqn:Cinv}
\vec{C}^{-1} =
\left(
\begin{array}{cc}
\L^{-\frac{1}{2}} & \vec{0} \\
\L^{-\frac{1}{2}}  \otimes \left( -\zeta^{-1} \M_\z^{\frac{1}{2}} \tilde{\z} \right) & \L^{-\frac{1}{2}}  \otimes  \zeta^{-1} \M_\z^{-\frac{1}{2}} 
\end{array}
\right) .
\end{align}
This factorization enables efficient manipulations in what follows.

\subsection{Bayesian inverse problem}
We formulate a Bayesian inverse problem to estimate $\d$ using the data pairs $\{\z_\ell,\y_\ell \}_{\ell=1}^N$. To facilitate the analysis we assume:
\begin{assumption}
$\{\z_\ell \}_{\ell=1}^N$ is a linearly independent set of vectors in $\R^n$.
\end{assumption}
This is usually satisfied in practice as one explores the model space in different directions since the number of model evaluations $N$ is small relative to the dimension of $\z$. For notational simplicity, we define
\begin{eqnarray}
\label{eqn:Aell}
\A_\ell = \begin{pmatrix} \I_{m} & \I_{m} \otimes \z_\ell^T \M_z \end{pmatrix} \in \R^{m \times p}, \qquad \ell=1,2,\dots,N,
\end{eqnarray}
so that $\d(\z_\ell,\t)=\A_\ell \t$, and the concatenation of these matrices
\begin{eqnarray*}
\A = \left(
\begin{array}{c}
\A_1 \\
\A_2 \\
\vdots \\
\A_N
\end{array}
\right) \in \R^{mN \times p} .
\end{eqnarray*}
$\A \t \in \R^{mN}$ corresponds to the evaluation of $\d(\z,\t)$ for all of the controller data $\{ \z_\ell \}_{\ell=1}^N$. In an analogous fashion, define $\vec{b} \in \R^{mN}$ by stacking $\y_\ell$, $\ell=1,2,\dots,N$, into a vector. Then we seek $\A \t \approx \vec{b}$.

To enable a closed form expression for the posterior we consider an additive Gaussian noise model with mean $\vec{0}$ and covariance $\alpha \vec{I}$, where $\alpha>0$ is specified by the user. Typically the noise covariance in a Bayesian inverse problem is defined using knowledge of the data collection process. In cases where data is collected from simulation we do not have a clear definition of noise. However, $\d$ is defined as a linear approximation of $\S(\z)-\tilde{\S}(\z)$, which is, in general, a nonlinear operator.  Hence the noise may be interpreted as the approximation error due to the linearization of the discrepancy. Nonetheless, one can show that the range space of $\A$ equals $\R^{mN}$, so there exists infinitely many $\t \in \R^p$ such that $\A \t=\vec{b}$. The choice of $\alpha$ dictates the weight given to the data misfit relative to the prior. Taking a small $\alpha$ will drive the inverse problem toward interpolation with little regard for the prior. We may experiment with different $\alpha$'s at a modest computational cost.

Given our formulation with Gaussian prior and noise models, and the linearity of $\d(\z,\t)$, the posterior is Gaussian with a negative log probability density (unnormalized) function 
\begin{eqnarray}
\label{eqn:neg_log_likelihood}
\frac{1}{2 \alpha} \left( \A \t - \vec{b} \right)^T \left( \A \t - \vec{b} \right) + \frac{1}{2} \t^T \M_\t \t .
\end{eqnarray}
The posterior mean and covariance are given by
\begin{eqnarray}
\label{eqn:post_mean_product}
\overline{\t} = \frac{1}{\alpha} \vec{\Sigma} \A^T \vec{b} \qquad \text{and} \qquad \vec{\Sigma} = \left( \M_\t + \frac{1}{\alpha} \A^T \A \right)^{-1},
\end{eqnarray}
respectively. Next we derive an explicit and computationally efficient expression for $\overline{\t}$ and its image under the post-optimality sensitivity operator.

\subsection{Expression for the posterior mean}
The posterior mean, given in~\eqref{eqn:post_mean_product}, involves the posterior covariance matrix which in turn involves the inverse of a sum of matrices in $\R^{p \times p}$. Since these matrices cannot be manipulated directly, we seek to derive a computable expression for the posterior covariance $\vec{\Sigma}$ and in turn the posterior mean $\overline{\t}$. Emphasis is given to preserving the Kronecker product structure throughout our analysis so that computations in $\R^p$ are avoided. The derivation is described in three steps in the following subsections: 
\begin{enumerate}
\item factorize $\A$ in a $\M_\t$ orthogonal basis to rewrite $\vec{\Sigma}^{-1}$,
\item invert $\vec{\Sigma}^{-1}$ by exploiting orthogonality to invert a sum,
\item compute $\overline{\t} = \frac{1}{\alpha} \vec{\Sigma} \A^T \vec{b}$ using Kronecker algebra.
\end{enumerate}

\subsubsection*{Factorize $\A$}
To rewrite the inverse posterior covariance, $ \M_\t + \frac{1}{\alpha} \A^T \A$, in a form amenable for inversion, we decompose $\A$ with respect to a chosen inner product. In particular, the generalized singular value decomposition (GSVD) of $\A$ is used with the inner product on the column space weighted by $\M_\t$. The GSVD is given by 
\begin{align*}
\A = \vec{\Xi} \vec{\Phi} \vec{\Psi}^T \M_\t
\end{align*}
where  $\vec{\Phi}$ is the diagonal matrix of singular values, and $\vec{\Xi}$ and $\vec{\Psi}$ are matrices containing the left and right singular vectors, respectively, which satisfy $\vec{\Xi}^T \vec{\Xi} = \I$ and $\vec{\Psi}^T \M_\t \vec{\Psi} = \I$. To determine the singular vectors, note that $\A \M_\t^{-1} \A^T=\vec{\Xi} \vec{\Phi}^2 \vec{\Xi}^T$. Using the expressions~\eqref{eqn:Aell} for $\A$ and~\eqref{eqn:Cinv} for $\M_\t^{-1}=\C^{-T} \C^{-1}$, we can write
\begin{eqnarray*}
\A \M_\t^{-1} \A^T =\vec{G} \otimes \L^{-1}
\end{eqnarray*}
where
\begin{align*}
& \G= \e \e^T + \zeta^{-2} (\vec{Z} - \tilde{\z} \e^T)^T \M_\z (\vec{Z}-\tilde{\z} \e^T) \in \R^{N \times N},
\end{align*}
$\vec{Z} = \begin{pmatrix} \z_1 & \z_2 & \dots & \z_N \end{pmatrix}$ is the matrix of controller data, and $\e \in \R^N$ is the vector of ones. Hence the left singular vectors $\vec{\Xi}$ correspond to the eigenvectors of $\vec{G} \otimes \L^{-1}$ and the squared singular values $\vec{\Phi}^2$ correspond to the eigenvalues.

Denoting the eigenvectors and eigenvalues of $\vec{G}$ and $\L$ with $(\vec{g}_i,\lambda_i)$ and $(\vec{l}_j,\rho_j)$, respectively, and recalling properties of the eigenvalue decomposition of a Kronecker product, we observe that the squared generalized singular values of $\vec{A}$ are given by $\frac{\lambda_i}{\rho_j}$ and are associated with the left singular vectors $\vec{\xi}_{i,j}=\vec{g}_i \otimes \vec{l}_j$. Rewriting the GSVD to solve for $\vec{\Psi}$, the right singular vector associated with $\vec{\xi}_{i,j}$ is given by 
 \begin{eqnarray}
 \label{eqn:psi_ij}
\vec{\psi}_{i,j}= \frac{1}{\sqrt{\lambda_i \rho_j}} 
\left(
 \begin{array}{cc}
s_i \vec{l}_j \\
 \vec{l}_j \otimes \zeta^{-2} \vec{w}_i
 \end{array}
 \right), \qquad i=1,2,\dots,N \ \ j=1,2,\dots,m,
 \end{eqnarray}
where
\begin{eqnarray}
\label{eqn:w_vecs}
\vec{w}_i = \vec{Z} \vec{g}_i - (\e^T \vec{g}_i) \tilde{\z} \qquad \text{and} \qquad s_i = (\e^T \vec{g}_i) - \zeta^{-2} \vec{w}_i^T \M_\z \tilde{\z} .
\end{eqnarray}

\subsubsection*{Invert $\vec{\Sigma}^{-1}$} 
Given this decomposition of $\A$ and our previously noted decomposition $\M_\t = \vec{C} \vec{C}^T$, we express $\vec{\Sigma}^{-1}$ as
\begin{eqnarray*}
\label{eqn:sigma_inv_factor}
\vec{\Sigma}^{-1} = \frac{1}{\alpha} \vec{C} \vec{X} \vec{C}^T 
\end{eqnarray*}
where
\begin{eqnarray}
\label{eqn:X}
\vec{X} =  \alpha \vec{I} + \vec{C}^T \vec{\Psi} \vec{\Phi}^2 \vec{\Psi}^T \vec{C} .
\end{eqnarray}

This factorization implies that $\vec{\Sigma}=\alpha \vec{C}^{-T} \vec{X}^{-1} \vec{C}^{-1}$ and facilitates easy manipulation because $\X$ is a combination of diagonal and orthogonal matrices. Applying the Sherman-Morrison-Woodbury formula to~\eqref{eqn:X}, we have
\begin{eqnarray*}
\vec{X}^{-1} = \frac{1}{\alpha} \left( \vec{I} - \vec{C}^T \vec{\Psi} \vec{D}\vec{\Psi}^T \vec{C} \right)
\end{eqnarray*}
where $\vec{D}\in \R^{mN \times mN}$ is a diagonal matrix whose entries are given by $\frac{\lambda_i}{\lambda_i + \alpha \rho_j}.$
Algebraic simplifications yield
\begin{eqnarray}
\label{eqn:Sigma}
\vec{\Sigma} = \M_\t^{-1} - \vec{\Psi} \vec{D}\vec{\Psi}^T .
\end{eqnarray}

The posterior covariance may be interpreted as the prior covariance $\M_\t^{-1}$ with its uncertainty reduced in the directions of the columns of $\vec{\Psi}$~\eqref{eqn:psi_ij} with weight $\frac{\lambda_i}{\lambda_i + \alpha \rho_j}$, a combination of the eigenvalues from $\L$ and $\G$, along with the noise covariance $\alpha$. 

\subsubsection*{Compute the posterior mean}
Recalling the form of the posterior mean in~\eqref{eqn:post_mean_product} and covariance in~\eqref{eqn:Sigma}, and computing the product
\begin{align*}
\A^T \b= \sum\limits_{\ell=1}^N 
\left(
\begin{array}{c}
\y_\ell \\
\y_\ell \otimes \M_\z \z_\ell
\end{array}
\right),
\end{align*}
 we have
\begin{align*}
\overline{\t} = \frac{1}{\alpha} \sum\limits_{\ell=1}^N \left( \M_\t^{-1} - \vec{\Psi} \vec{D}\vec{\Psi}^T \right)
\left(
\begin{array}{c}
\y_\ell \\
\y_\ell \otimes \M_\z \z_\ell
\end{array}
\right) .
\end{align*}

Using $\M_\t^{-1} = \C^{-T} \C^{-1}$, recalling the form of $\vec{\Psi}$'s columns~\eqref{eqn:psi_ij}, and writing matrix-vector products with $\vec{\Psi} \vec{D}\vec{\Psi}^T$ in terms of linear solves involving $\L$, we arrive at the posterior mean expression
\begin{eqnarray*}
\overline{\t} = \frac{1}{\alpha} \sum\limits_{\ell=1}^N 
\left[
\left(
\begin{array}{c}
a_\ell \u_\ell \\
 \u_\ell  \otimes  \zeta^{-2} (\z_\ell-\tilde{\z}) 
 \end{array}
\right)
-
\sum\limits_{i=1}^N b_{i,\ell}
\left(
\begin{array}{c}
s_i \vec{u}_{i,\ell} \\
 \vec{u}_{i,\ell} \otimes   \zeta^{-2} \vec{w}_i 
\end{array}
\right)
\right]
\end{eqnarray*}
involving the constants
\begin{align*}
 a_{\ell} = 1- \zeta^{-2}  \tilde{\z}^T\M_\z (\z_\ell-\tilde{\z}) \qquad \text{and} \qquad  b_{i,\ell} = \zeta^{-2} (\z_\ell-\tilde{\z})^T \M_\z \vec{Z} \vec{g}_i + (\e^T \vec{g}_i) a_\ell
\end{align*}
and the vectors
\begin{eqnarray*}
\vec{u}_\ell = \L^{-1} \y_\ell \qquad \text{and} \qquad \vec{u}_{i,\ell} = \left( \alpha \L + \lambda_i\I \right)^{-1} \u_\ell .
\end{eqnarray*}
We observe that $\overline{\t}$ depends on linear solves \footnote{In Section~\ref{sec:numerics}, we demonstrate how to efficiently invert $\L$ and $ \left( \alpha \L + \lambda_i\I \right)$ using a generalized singular value decomposition when $\L$ is defined as the square of an elliptic operator.} involving $\L$ and $\alpha \L + \lambda_i\I$ where the right hand sides arise from the data $\{\y_\ell\}_{\ell=1}^N$. The posterior mean can be interpreted as a linear combination of the discrepancy data preconditioned by the prior weighted by the noise covariance $\alpha$ and the eigenpairs $(\lambda_i,\g_i)$ from the controller data informed matrix $\G$. 

\subsection{Optimal solution update}
\label{subsec:post_model_discrepancy}
To propagate $\overline{\t}$ through the post-optimality sensitivity operator $\nabla_{\t} \F(\vec{0})$, we note the Kronecker product structure of both $\B$ and $\overline{\t}$ which allows for computation of the post-optimality sensitivities using solves in $\R^m$ and $\R^n$. The updated optimal solution is given by $\tilde{\z} - \H^{-1} \B \overline{\t}$ where
\begin{align}
\label{eqn:B_thetabar}
\B \overline{\t} &=  \frac{1}{\alpha} \nabla_{\z} \tilde{\S}^T \nabla_{\u, \u} \J \left[ \sum\limits_{\ell=1}^N \left(  \u_\ell-  \sum\limits_{i=1}^N b_{i,\ell} (\e^T \vec{g}_i) \vec{u}_{i,\ell} \right)\right] \\
& + \frac{1}{\alpha} \sum\limits_{\ell=1}^N ( \nabla_\u J \u_\ell) \zeta^{-2} \M_\z ( \z_\ell-\tilde{\z}) \nonumber \\
&- \frac{1}{\alpha} \sum\limits_{\ell=1}^N \sum\limits_{i=1}^N b_{i,\ell} (\nabla_\u \J \vec{u}_{i,\ell}) \zeta^{-2} \M_\z \vec{w}_i \nonumber .
\end{align}
Equation~\eqref{eqn:B_thetabar} provides a systematic and interpretable combination of high-fidelity data preconditioned by the prior in $\vec{u}_\ell$ and $\vec{u}_{i,\ell}$, the low-fidelity model in $\tilde{\vec{S}}_\z^T$, and the optimization objective in $\nabla_\u \J$ and $\nabla_{\u,\u} \J$. 

\section{Specifying the prior} \label{sec:prior_discrepancy}
The proposed approach requires specification of the matrix $\L$, the scalar $\zeta$, and the scalar $\alpha$. This is necessary because the inverse problem to estimate $\d$ is underdetermined and hence is posed in a Bayesian formulation. These parameters play an important role imposing domain expertise to constraint the discrepancy. This section explores a strategy to choose them based on physical characteristics of the problem.

The noise covariance $\alpha>0$ normalizes the data misfit in the likelihood thus weighting the data fit relative to the prior~\eqref{eqn:neg_log_likelihood}. Typically in Bayesian inverse problems $\alpha$ is specified based on knowledge of the data fidelity. The noise covariance in our framework corresponds to the approximation error due to linearization of the discrepancy. As a rule of thumb, $\alpha$ may be chosen based on the magnitude of the discrepancy data. For instance, if the discrepancy data $\y_\ell$ has magnitude $\mathcal O(1)$ then $\alpha=10^{-2}$ normalizes the difference $\y_\ell-\d(\z_\ell,\t)$ by $\sqrt{\alpha}=10^{-1}$ or equivalently $10\%$ of the discrepancy magnitude. Comparing $\d(\z_\ell,\overline{\t})$ with $\y_\ell$ provides an easy check that the choice was reasonable.

Sampling the prior discrepancy aids in determining $\L$ and $\zeta$ by analyzing the physical characteristics of the discrepancy samples. The matrix $\L$ imposes spatiotemporal characteristics (variance, correlation length, boundary conditions, etc.) on the state space and the scalar $\zeta$ defines length scales on the controller space to affect how the discrepancy varies with respect to $\z$. A sampling procedure is detailed below to facilitate the specification of $\L$ and $\zeta$ independently of one another. 

Using the factorization $\M_\t^{-1} = \vec{C}^{-T} \vec{C}^{-1}$, a prior sample of $\vec{\Theta}$ is given by $\vec{C}^{-T} \vec{\omega}_p$, where $\vec{\omega}_p$ is a $p$ dimensional standard normal random vector. Evaluating the model discrepancy at such samples yields
\begin{align*}
\d(\z,\vec{C}^{-T} \vec{\omega}_p) = \L^{-\frac{1}{2}} \vec{\omega}_0 + \zeta^{-1} \L^{-\frac{1}{2}} \vec{\Omega} \M_\z^{\frac{1}{2}} (\z-\tilde{\z})
\end{align*}
where $\vec{\omega}_0 \in \R^m$ and $\vec{\Omega} \in \R^{m \times n}$ have entries that are independent identically distributed (i.i.d) samples from a standard normal distribution. 

To avoid forming the dense matrix $\vec{\Omega} \in \R^{m \times n}$, we evaluate $\d(\z,\vec{C}^{-T} \vec{\omega})$ at particular $\z$ values and visualize the resulting discrepancy, which is an element of the state space. Let $\z_r \in \R^n$ denote a ``reference point" corresponding to a physically reasonable controller. We will evaluate the discrepancy at a sequence of points $\hat{\z}_k=\tilde{\z} + \frac{k}{K}(\z_r-\tilde{\z})$, $k=0,1,\dots,K$, where $K$ is an integer. Then for a fixed $\vec{\Theta}$ sample, or equivalently a fixed $\vec{\omega}_p$, we have
\begin{align*}
\d(\hat{\z}_k,\vec{C}^{-T} \vec{\omega}_p) = \L^{-\frac{1}{2}} \vec{\omega}_0 + \frac{k}{K} \zeta^{-1} \L^{-\frac{1}{2}} \vec{\Omega} \M_\z^{\frac{1}{2}} (\z_r-\tilde{\z}).
\end{align*}

Rather than forming the large matrix $\vec{\Omega} \in \R^{m \times n}$ of i.i.d. random numbers, observe that $ \vec{\Omega} \M_\z^{\frac{1}{2}} (\z_r-\tilde{\z})$ is a Gaussian random vector in $\R^m$ whose entires are i.i.d. Gaussians with mean $0$ and variance $(\z_r-\tilde{\z})^T \M_\z (\z_r-\tilde{\z})$. Thus to compute $\d(\hat{\z}_k,\vec{C}^{-T} \vec{\omega}_p)$ for each $k=0,1,\dots,K$, it is sufficient to apply $\L^{-\frac{1}{2}}$  to two random vectors $\vec{\omega}_0,\vec{\omega}_1 \in \R^m$, sampled with i.i.d. standard normal entries, and evaluate
\begin{align*}
\d(\hat{\z}_k,\vec{C}^{-T} \vec{\omega}_p) = \L^{-\frac{1}{2}} \vec{\omega}_0 + \frac{k}{K} \zeta^{-1} \sqrt{(\z_r-\tilde{\z})^T \M_\z (\z_r-\tilde{\z})}  \L^{-\frac{1}{2}} \vec{\omega}_1 .
\end{align*}
Our approach is summarized in Algorithm~\ref{alg:prior_samples} where $S$ samples are generated by looping over different standard normal random vectors $\vec{\omega}_0^s,\vec{\omega}_1^s$, $s=1,2,\dots,S$. Examination of $ \{ \L^{-\frac{1}{2}} \vec{\omega}_0^s \}_{s=1}^S $ and $\{ \L^{-\frac{1}{2}} \vec{\omega}_1^s \}_{s=1}^S$ aids in tuning of parameters in $\L$ so that the discrepancy realizations are physically reasonable. Similarly, $\zeta$ may be tuned based on variability in $\{\d(\hat{\z}_k, \vec{C}^{-T} \vec{\omega}_p^s) \}_{k=1}^K$ for each fixed $\vec{\omega}_p^s$, where $\vec{\omega}_p^s \in \R^p$ is not computed but is used to denote the sample corresponding to $\vec{\omega}_0^s, \vec{\omega}_1^s \in \R^m$.
\begin{algorithm}[ht!!]
	\caption{Compute samples of the prior discrepancy.}
\begin{algorithmic} [1] \label{alg:prior_samples}
		\STATE Input: $S \in \mathbb{N}, K \in \mathbb{N}$, and $\z_r \in \R^n$ 
		\STATE Compute $c=\sqrt{(\z_r-\tilde{\z})^T \M_\z (\z_r-\tilde{\z})} $
		\FOR{$s = 1,2,\dots,S$}
		\STATE Sample $\vec{\omega}_0,\vec{\omega}_1 \in \R^m$ with i.i.d. standard normal entries
		\STATE Compute $\L^{-\frac{1}{2}} \vec{\omega}_0$ and $\L^{-\frac{1}{2}} \vec{\omega}_1$
		\FOR{$k=0,1,\dots,K$}
		\STATE Evaluate $\d(\hat{\z}_k,\vec{C}^{-T} \vec{\omega}^s)=\L^{-\frac{1}{2}} \vec{\omega}_0 + \frac{k}{K} \zeta^{-1} c \L^{-\frac{1}{2}} \vec{\omega}_1$
		\ENDFOR
		\ENDFOR
	\end{algorithmic}
\end{algorithm}

\section{Numerical results} \label{sec:numerics}
We demonstrate our approach on an illustrative 1D control example, where model discrepancy arises from a failure to include advection, and on an application of fluid flow in 2D, where Stokes is used as the low-fidelity approximation of  Naiver-Stokes.

\subsection{Illustrative example}
A 1D steady state control problem is formulated that seeks an optimal source to achieve a target state profile. The problem is formulated as
\begin{align*}
 \min_{z} \frac{1}{2} \int_0^1 (\tilde{S}(z)(x)-T(x))^2 dx + \frac{\beta}{2} \int_0^1 z(x)^2 dx  
\end{align*}
where $T(x) = 50 - 30(x-0.5)^2$ is the target, $z:[0,1] \to \R$ is a source controller, $\beta = 10$ is a regularization coefficient, and $\tilde{S}(z)$ is the solution operator for the diffusion equation
\begin{align*}
& - \kappa u'' = z \qquad & \text{on } (0,1) \\
& \kappa u' = h u     \quad & \text{on }     \{0\} \\ 
 & -\kappa u' = h u  \quad & \text{on } \{1\} 
\end{align*}
with Robin boundary condition. The state and control spaces $\U$ and $\Z$ are $L^2(0,1)$. Our assumption is that the diffusion model is missing physics, specifically it lacks advection.  We consider the high-fidelity model to be
\begin{align*}
& - \kappa u'' + v u' = z \qquad & \text{on } (0,1) \\
& \kappa u' = h u     \quad & \text{on }     \{0\} \\ 
 & -\kappa u' = h u  \quad & \text{on } \{1\} 
\end{align*}
represented by solution operator $S(z)$. We take parameter values $\kappa=1$, $v=0.5$, and $h= 2$.

\subsubsection*{Optimization}
After discretizing with linear finite elements on a uniform mesh of $200$ nodes, i.e. $m=n=200$, we solve both high and low-fidelity optimization problems to illustrate the effect of model discrepancy,  which is attributed to the missing advection term. We observe a considerable difference in the optimal sources, shown in Figure~\ref{fig:heat_example_posterior_source_samples}. Since advection moves the state from left to right, including advection in the model results in a source skewed to the left.

\subsubsection*{Prior model discrepancy}

We define $\L = \vec{D} \vec{M}^{-1} \vec{D}$, where $\vec{D}$ is the discretization of an elliptic operator $\gamma (-\epsilon \Delta + \mathcal I)$, equipped with zero Neumann boundary conditions, and $\vec{M}$ is the mass matrix for the state discretization. The coefficients $\gamma$ and $\epsilon$ determine the magnitude and smoothness, respectively, of the prior discrepancy realizations. The use of an elliptic operator ensures interpretability of the hyper-parameters and computational convenience of efficient linear solves and matrix factorizations for sampling~\citep{ghattas_infinite_dim_bayes_1,ghattas_infinite_dim_bayes_2,stuart_inv_prob}. 

\begin{figure}[h]
\centering
  \includegraphics[width=0.3\textwidth]{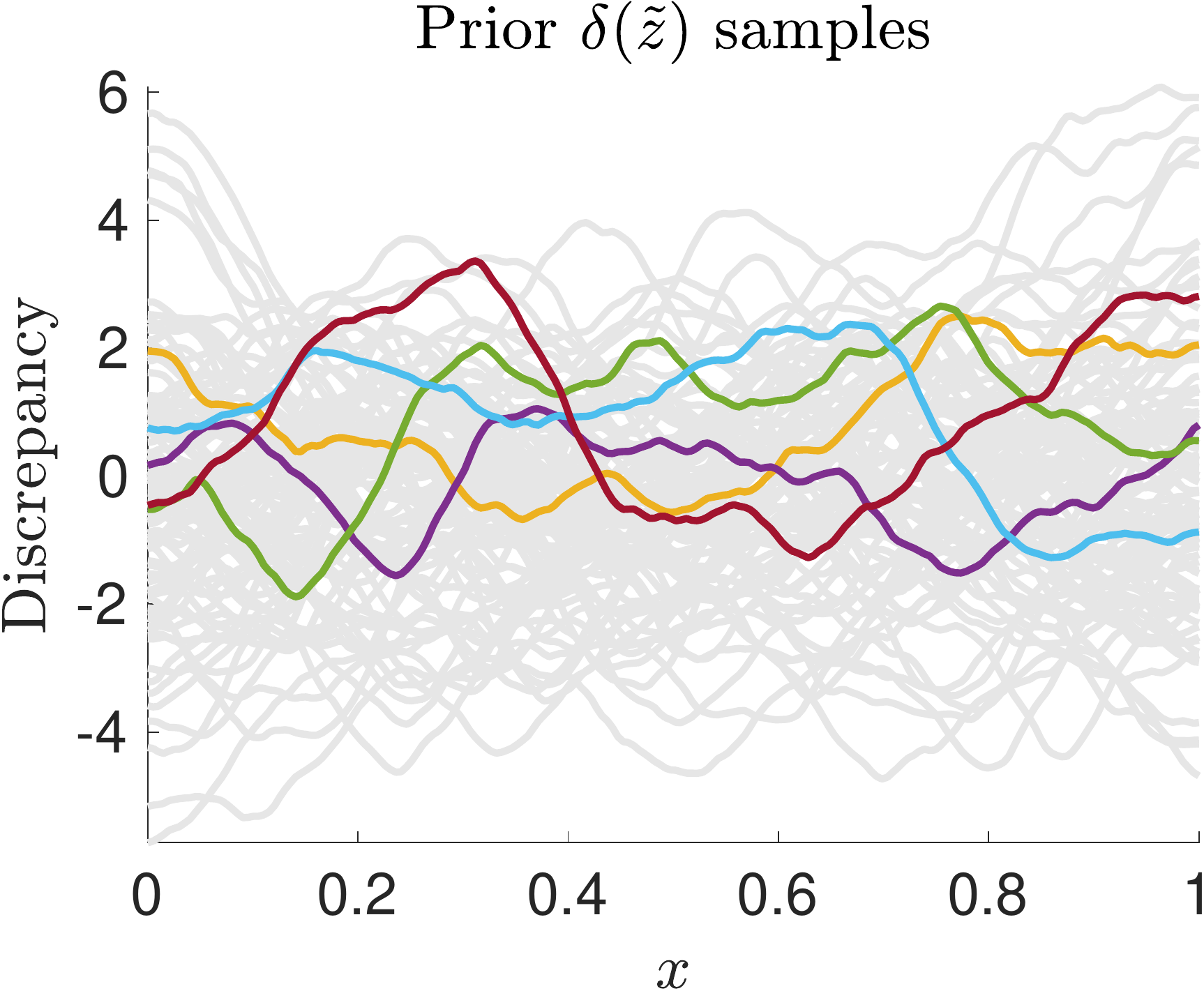}
    \includegraphics[width=0.3\textwidth]{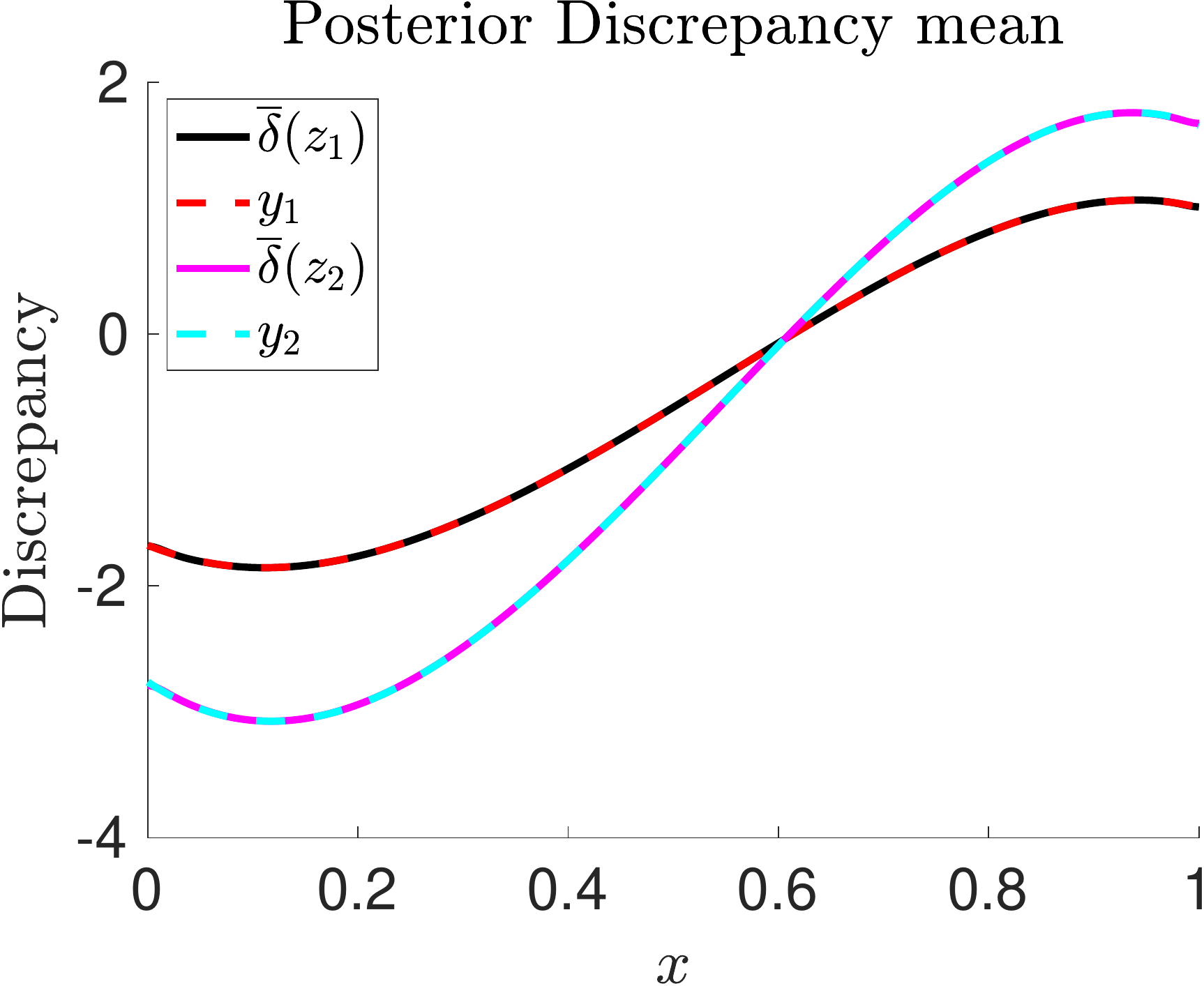}
  \caption{Left: 100 prior discrepancy samples (5 are colored for ease of visualization) evaluated at $\tilde{\z}$; right: mean posterior discrepancy evaluated at the two datapoints $\{ \z_1,\z_2\}$ compared with the discrepancy data $\{\y_1,\y_2\}$.}
  \label{fig:heat_example_prior_samples_1}
\end{figure}

Following Algorithm~\ref{alg:prior_samples}, we compute samples of the prior discrepancy to tune $\gamma$ and $\epsilon$ so that the samples are physically consistent with the magnitudes and smoothness inherent to this problem. The left panel of Figure~\ref{fig:heat_example_prior_samples_1} displays the samples corresponding to our choice of $\gamma =1$ and $\epsilon=10^{-2}$. The observed data $\y_\ell=\S(\z_\ell)-\tilde{\S}(\z_\ell)$ is shown in the right panel of Figure~\ref{fig:heat_example_prior_samples_1} to assess the discrepancy's magnitude, smoothness, and dependence on $\z$. 

\begin{figure}[h]
\centering
  \includegraphics[width=0.24\textwidth]{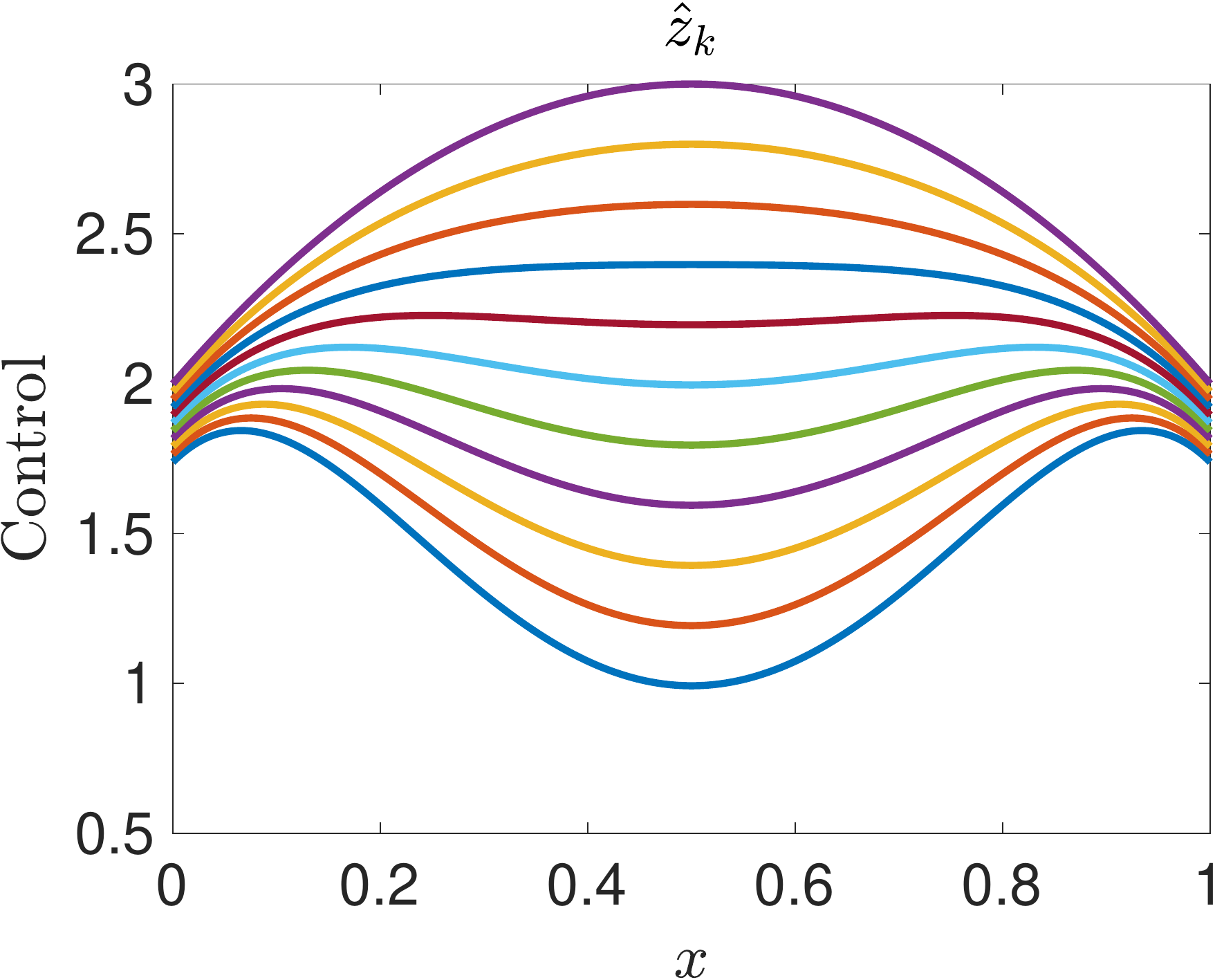}
   \includegraphics[width=0.24\textwidth]{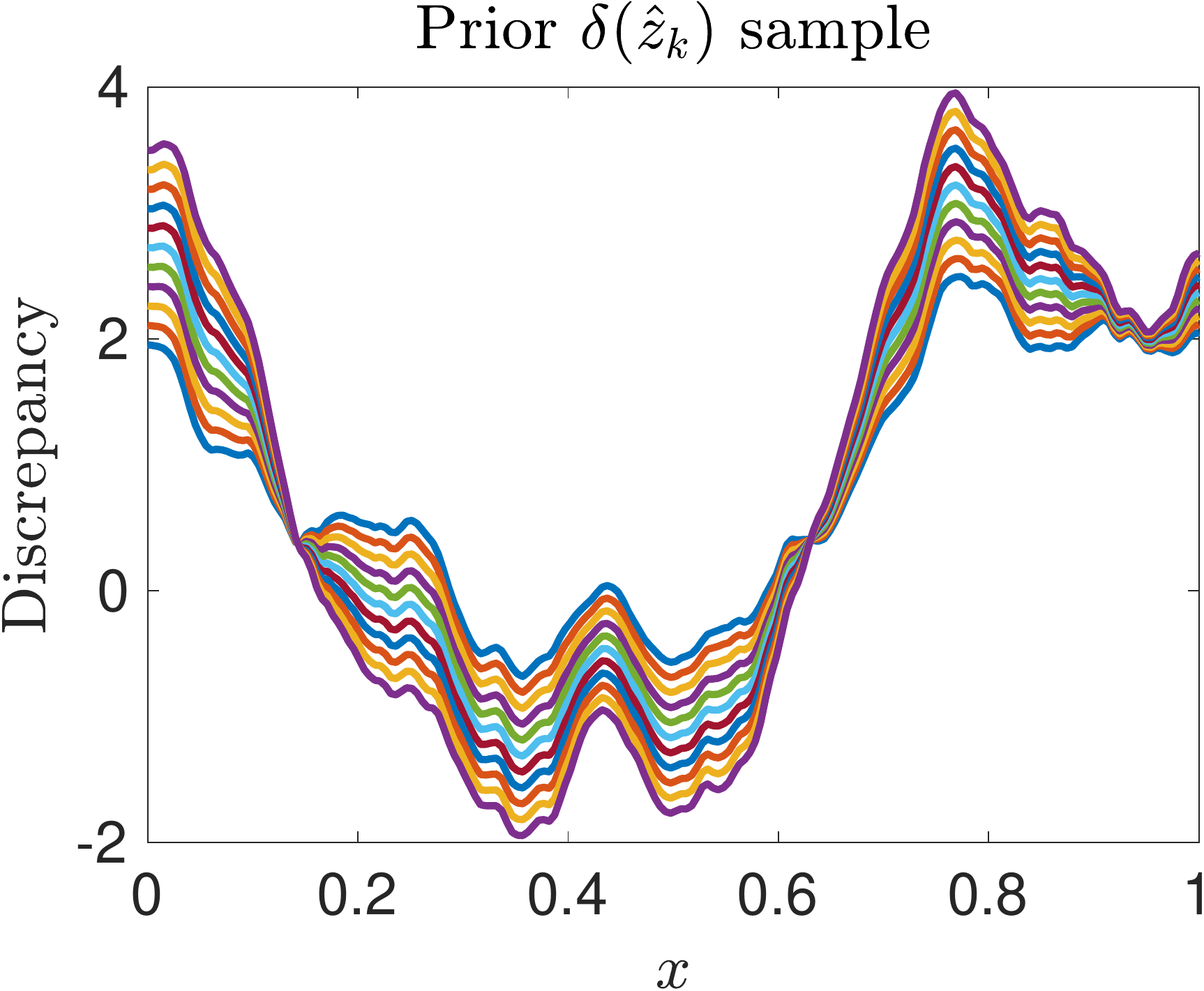}
    \includegraphics[width=0.24\textwidth]{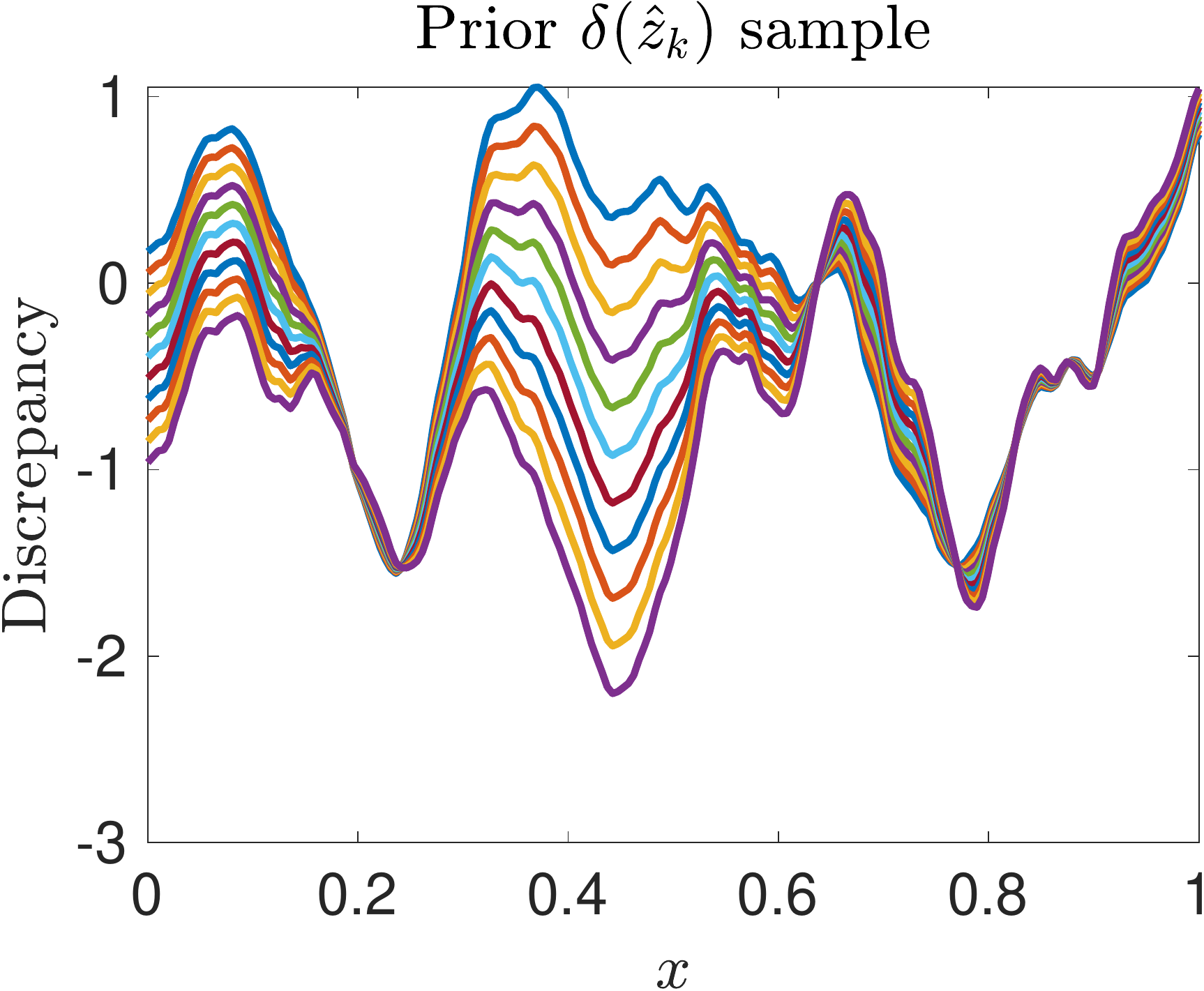}
     \includegraphics[width=0.24\textwidth]{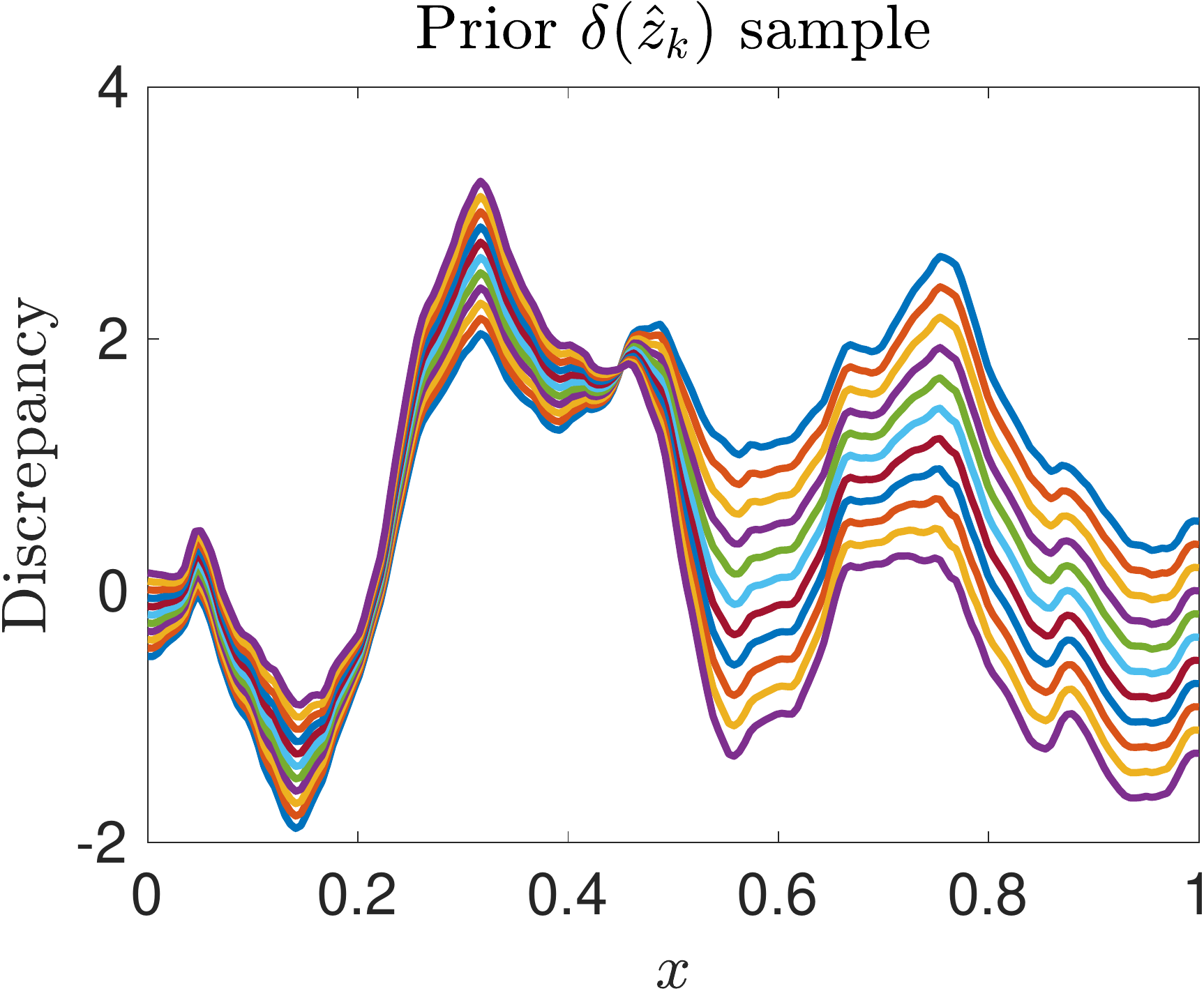}
  \caption{Leftmost panel: sequence of controllers $\hat{\z}_k=\tilde{\z} + \frac{k}{11}(\z_r-\tilde{\z})$, $k=0,1,\dots,11$; three rightmost panels: three prior discrepancy samples evaluated at $\hat{\z}_k$, $k=0,1,\dots,11$.}
  \label{fig:heat_example_prior_discrepancy_vary_z_samples}
\end{figure}

The control length  scale parameter $\zeta$ is determined by analyzing variability of the prior samples with respect to $\z$, as described in Section~\ref{sec:prior_discrepancy}. Figure~\ref{fig:heat_example_prior_discrepancy_vary_z_samples} shows $\{\hat{\z}_k\}_{k=0}^{11}$ in the leftmost panel and three samples $\{ \d(\hat{\z}_k,\vec{C}^{-T} \vec{\omega}_p^s)\}_{k=0}^{11}$, $s=1,2,3$, in the three rightmost panels. The scalar $\zeta$ may be tuned according to the expected rate of variation of the discrepancy with respect to $\z$. In our case, $\zeta=2$. 

\subsubsection*{Posterior model discrepancy}
We evaluate the high-fidelity model twice to generate $N=2$ data pairs $\{\z_\ell,\y_\ell\}_{\ell=1}^2$ and take a noise variance $\alpha=0.01$. The right panel of Figure~\ref{fig:heat_example_prior_samples_1} displays the mean posterior discrepancy evaluated at $\z_1$ and $\z_2$. We observe that it matches the data well. This is unsurprising given that $\S(\z)-\tilde{\S}(\z)$ is a linear function of $\z$ and $\alpha$ was taken small enough to encourage trust in the data. 

\subsubsection*{Optimal solution update}
\begin{figure}[h]
\centering
  \includegraphics[width=0.4\textwidth]{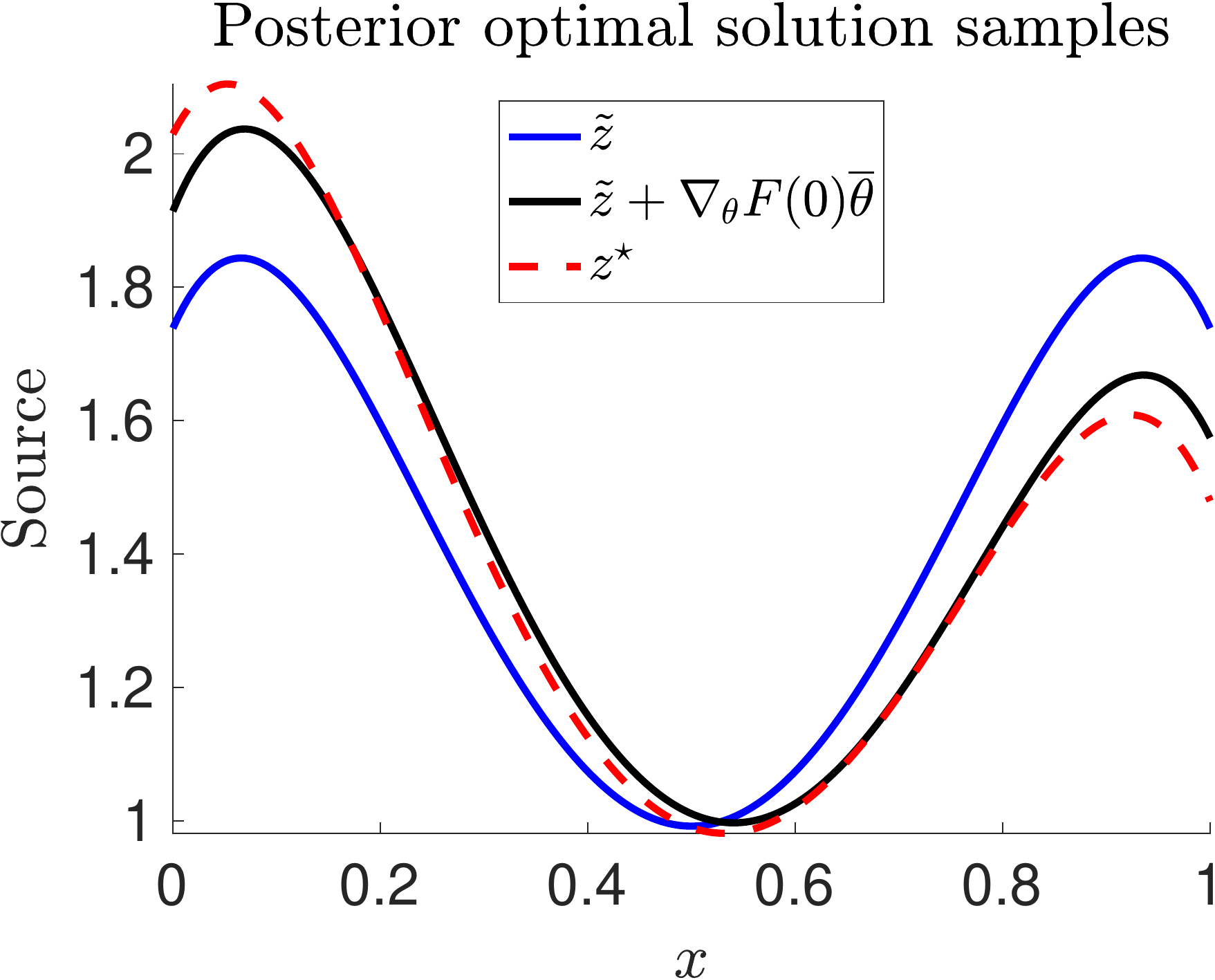}
  \caption{The low and high-fidelity solutions $\tilde{\z}$ and $\z^\star$ are shown by the blue and broken red lines, respectively. The updated optimal solution is given by the black line.}
  \label{fig:heat_example_posterior_source_samples}
\end{figure}

Propagating the posterior model discrepancy mean through the post-optimality sensitivity operator, we arrive at the updated optimal solution shown in Figure~\ref{fig:heat_example_posterior_source_samples}. We observe that the updated optimal solution is much closer to the high-fidelity solution $\z^\star$ than the low-fidelity solution $\tilde{\z}$. Hence using the low-fidelity optimization problem along with two high-fidelity forward solves, we are able to find a good approximation of the high-fidelity optimal solution. 

\subsection{Control of viscous fluid flow}
In this subsection our approach is demonstrated on control of viscous fluid flow. In particular, we consider the high-fidelity model to be governed by the Naiver-Stokes equation and the low-fidelity model represented by the Stokes equation. The goal is find the optimal distributed controller\footnote{We use the notation $\vec{z}$ to denote the controller before discretization since it is vector-valued. This should not be confused with the controller coordinates $\vec{z} \in \R^n$ used throughout the article.} $\vec{z}=(z_x,z_y)$, defined on a subset of the domain, that minimizes the vertical fluid flow. The low-fidelity optimization problem is
\begin{align*}
 \min_{z} \frac{1}{2} \int_{\chi} \tilde{v}_y(z)^2  + \frac{\beta}{2} \int_{\Omega_z} (z_x^2 + z_y^2) 
\end{align*}
where $\tilde{\vec{v}}(z)=(\tilde{v}_x(z),\tilde{v}_y(z))$ is the velocity solution operator for the Stokes equation
\begin{align*}
& -\mu \nabla \vec{v} + \nabla p = \vec{g} + \vec{z} & \text{on } \Omega \\
& \nabla \cdot \vec{v} = 0 & \text{on } \Omega 
\end{align*}
where $\Omega =(0,1)^2$, $\Omega_z=(0.2,.8)\times (0.1,0.4)$, and $\chi = (0,1) \times (0,0.5)$. The controller is defined to be zero on $\Omega \setminus \Omega_z$. The flow is driven by the effect of gravity $\vec{g}=(0,9.81)$, viscosity $\mu=0.5$, and the inflow boundary conditions 
\begin{eqnarray*}
v_x(0,y) = 6y(1-y) \qquad \text{and} \qquad v_y(x,1) = -2 \sin(2 \pi x)^2.
\end{eqnarray*}
The state and control spaces $\U$ and $\Z$ are $L^2(\Omega)$ and $L^2(\Omega_z)$, respectively. We solve the optimization problem with regularization coefficient $\beta=10^{-4}$. The high-fidelity Naiver-Stokes model
\begin{align*}
& -\mu \nabla \vec{v} + (\vec{v} \cdot \nabla) \vec{v} + \nabla p = \vec{g} + \vec{z} & \text{on } \Omega \\
& \nabla \cdot \vec{v} = 0 & \text{on } \Omega 
\end{align*}
includes the nonlinear convective term omitted in the Stokes equation.

\subsubsection*{Optimization}

\begin{figure}[h]
\centering
\includegraphics[width=0.32\textwidth]{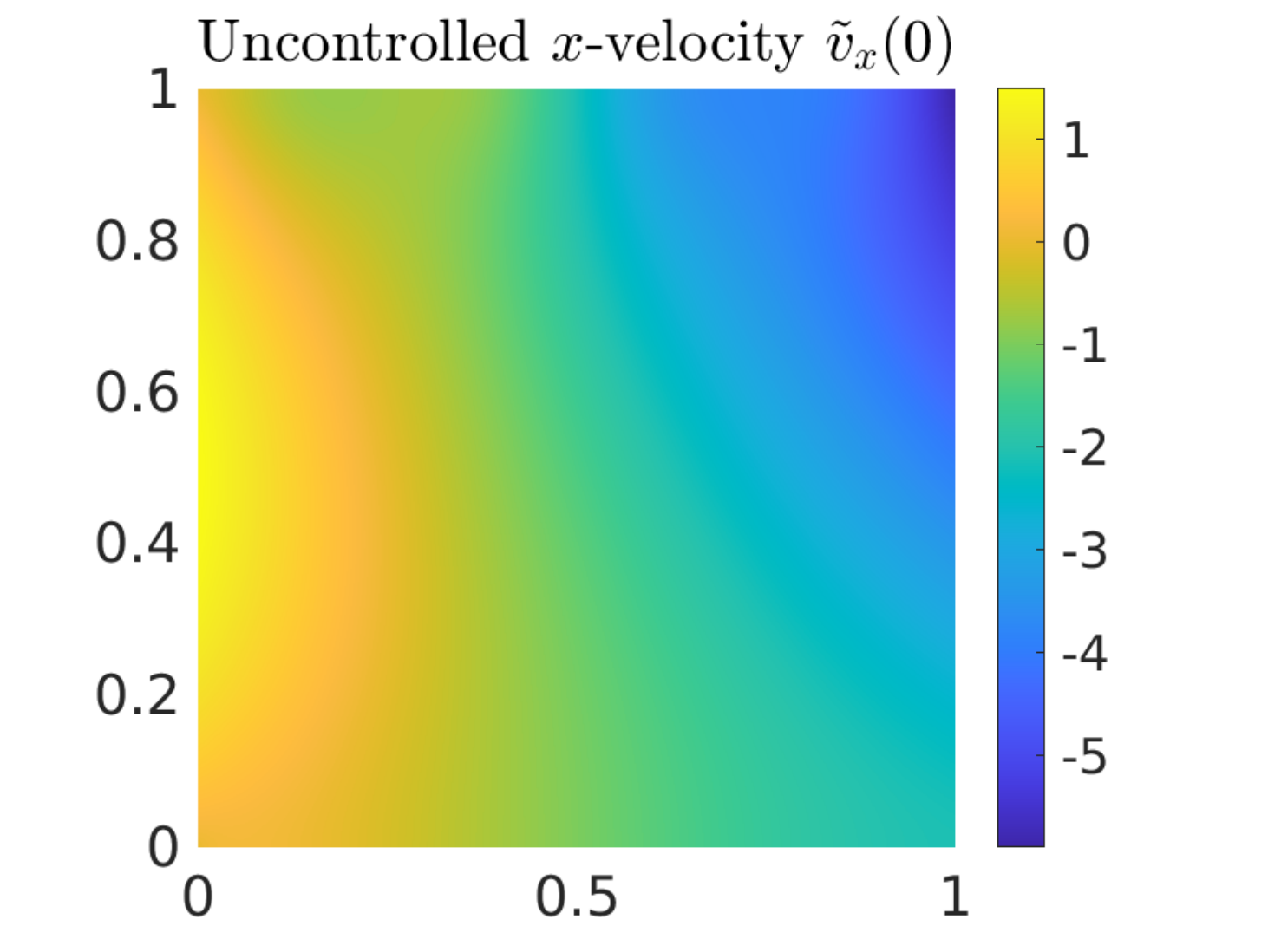}
\includegraphics[width=0.32\textwidth]{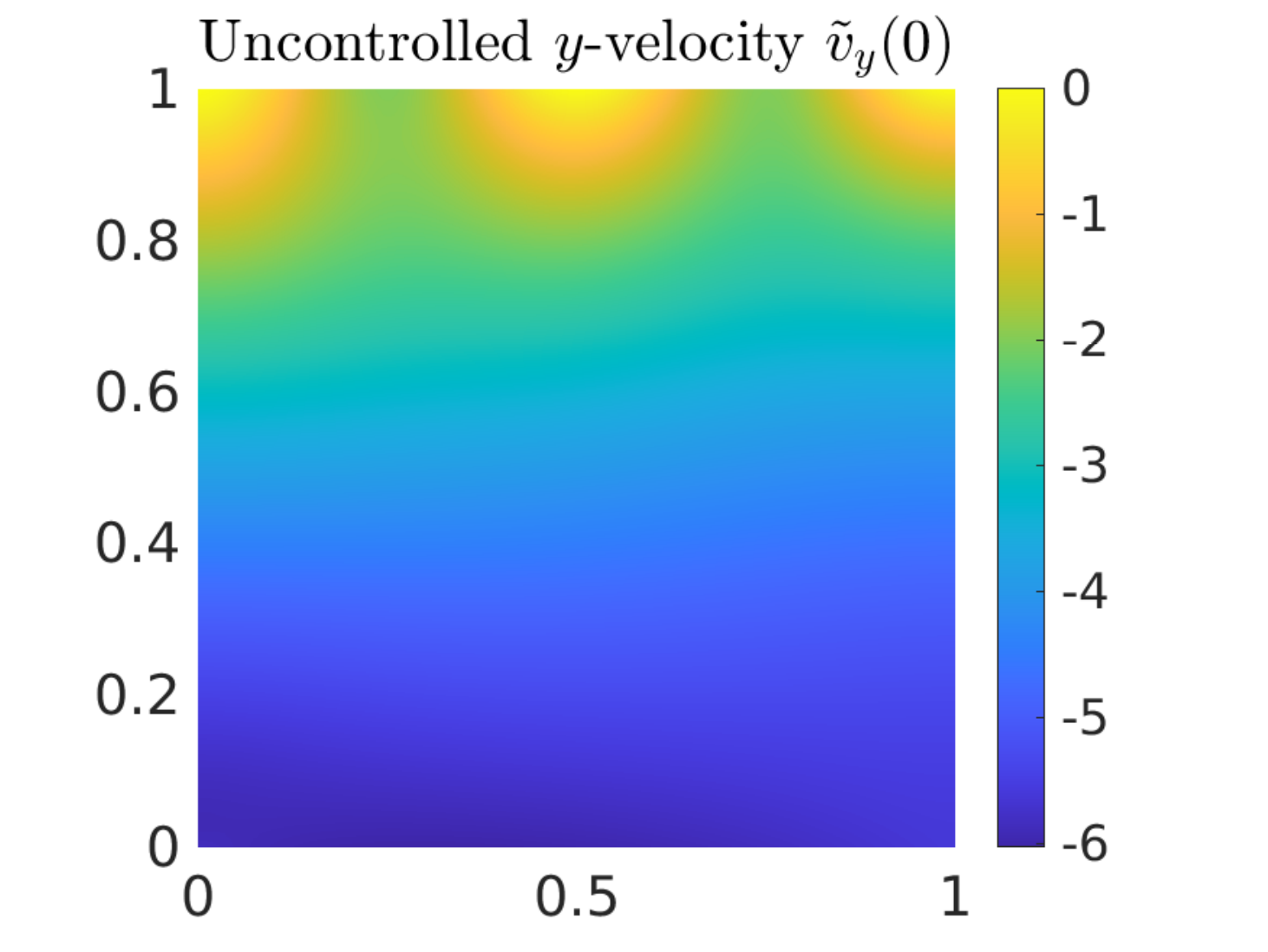}
 \includegraphics[width=0.32\textwidth]{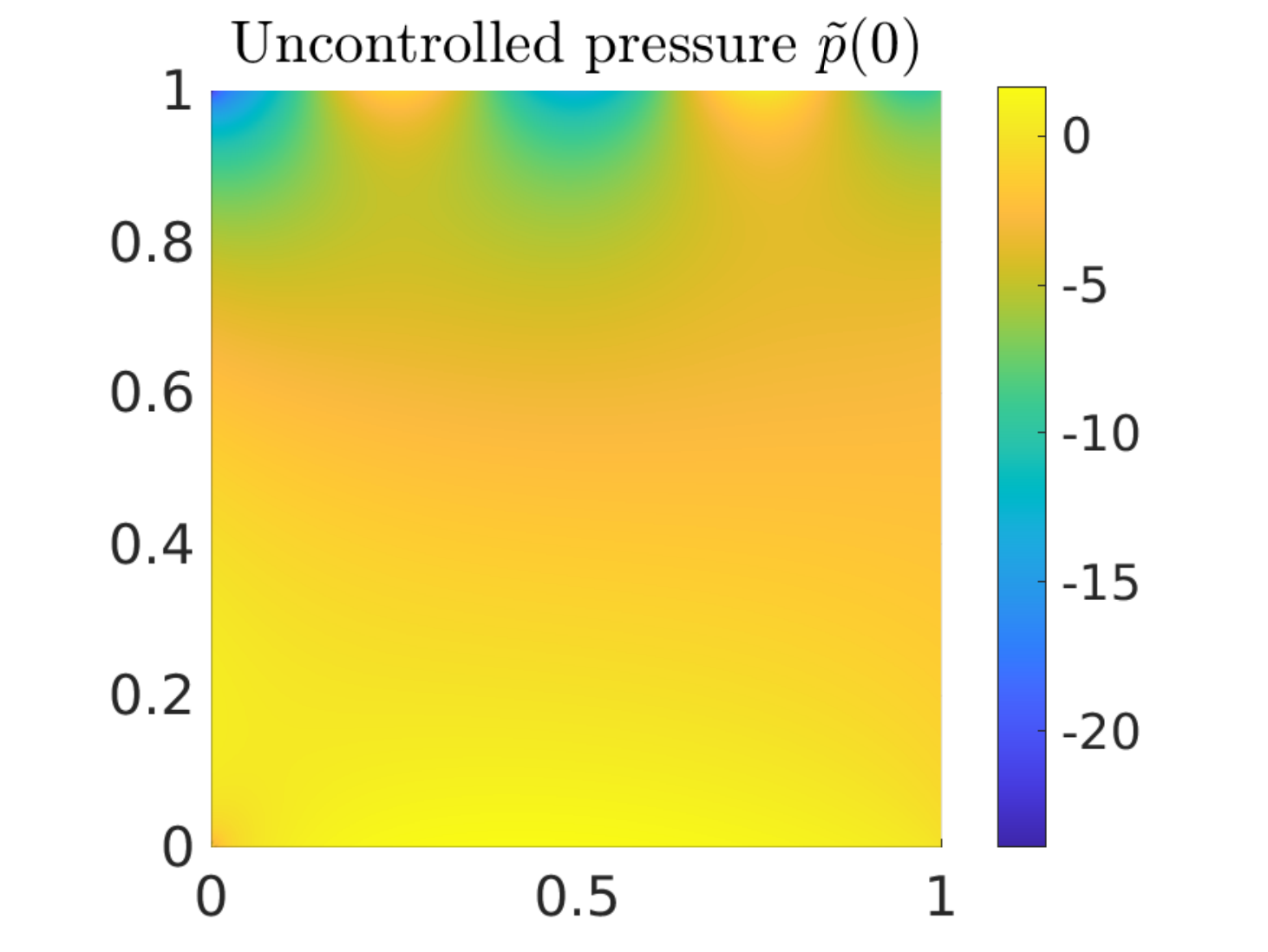} \\
\includegraphics[width=0.32\textwidth]{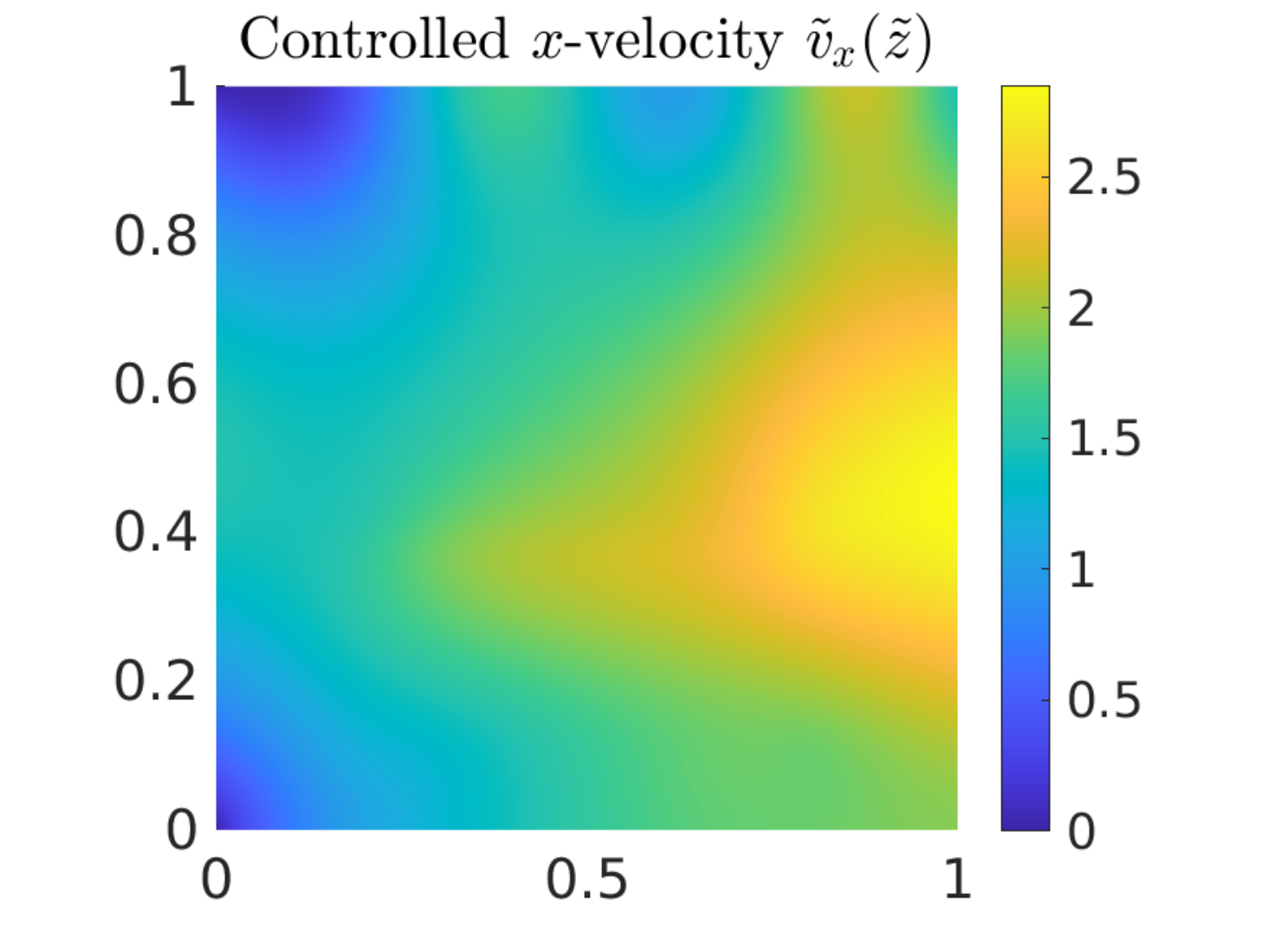}
\includegraphics[width=0.32\textwidth]{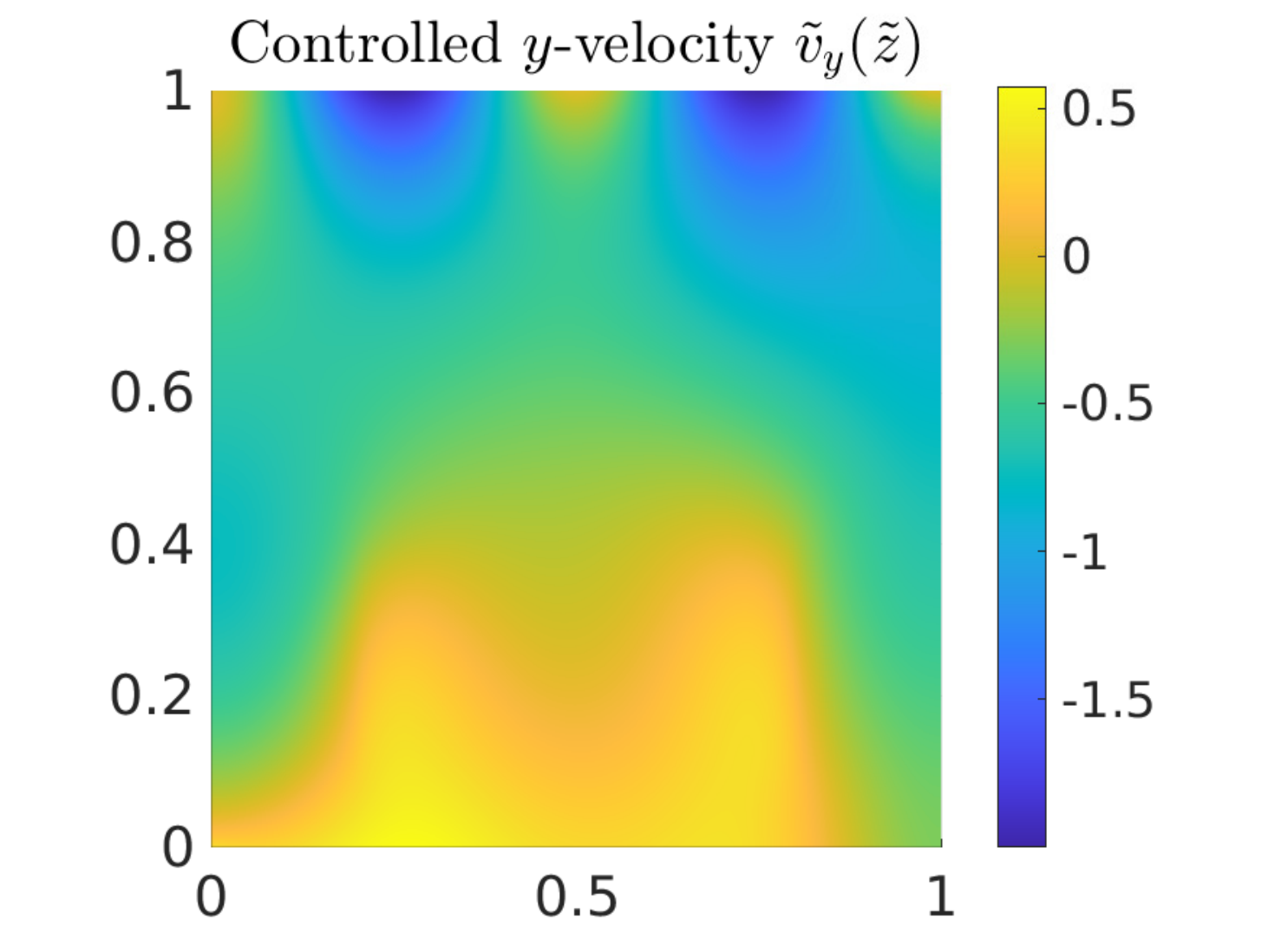}
 \includegraphics[width=0.32\textwidth]{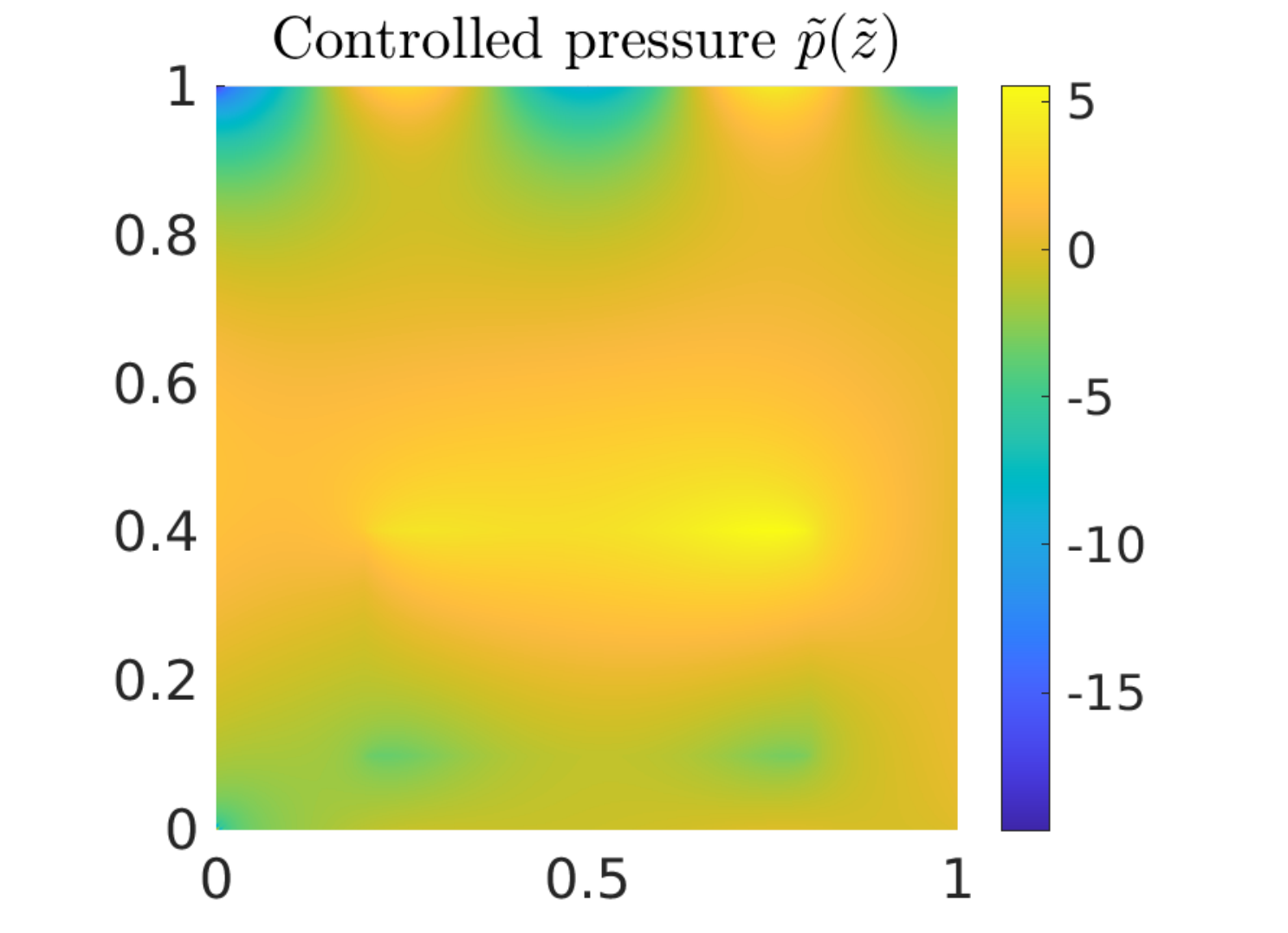}
  \caption{Uncontrolled (top) and controlled (bottom) state solution for the Stokes equation. The $x$-velocity $v_x$, $y$-velocity $v_y$, and pressure $p$ are plotted from left to right. The uncontrolled $y$-velocity has a magnitude around $6$ everywhere in $\chi$ whereas the controlled $y$-velocity has a magnitude close to $0$.}
  \label{fig:stokes_example_state_solution}
\end{figure}

The three state variables $v_x$, $v_y$, and $p$ are discretized with Taylor-Hood elements yielding $m=365622$ degrees of freedom. Figure~\ref{fig:stokes_example_state_solution} displays the uncontrolled and optimally controlled solutions of the Stokes equation. The inflow from the left and above, combined with the force of gravity, produce a flow field with high downward velocity in the lower region of the domain. The speed of this flow is minimized and the efficacy of the controller (for the Stokes equation) is apparent from the center column of Figure~\ref{fig:stokes_example_state_solution}. However, the Naiver-Stokes solution evaluated at the optimal controller exhibits faster flows.

\subsubsection*{Prior model discrepancy} 
Similar to the previous example, we define the prior covariance $\L^{-1}$ as the square of an inverse elliptic differential operator. However, to mitigate the effect of high boundary variance in the prior we follow~\citep{stadler_cov_bc} by imposing a Robin boundary condition which is optimized to reduce the boundary effect. Because this condition is not sufficient to completely eliminate the inflated variance near the boundary, a rescaling of the covariance based on a point-wise variance estimate using $100$ samples is also implemented.

The covariance is $\L^{-1}=\gamma^2(\vec{E}^{-\frac{1}{2}} \vec{D} \vec{M}^{-1} \vec{D} \vec{E}^{-\frac{1}{2}})^{-1}$, where $\vec{D}$ is the discretization of the elliptic operator $(-\epsilon \Delta + \mathcal I)$, $\gamma=2$ and $\epsilon=0.2$, and $\vec{E}$ is a diagonal matrix that approximates the diagonal of $( \vec{D} \vec{M}^{-1} \vec{D})^{-1}$, and $\vec{M}$ is the mass matrix for the state discretization. We omit showing prior samples for conciseness, but emphasize that the interpretably of the elliptic operator facilitates specification of $\gamma$ and $\epsilon$ to ensure that the prior is reflective of the physical characteristics of the discrepancy. We use $N=1$ high-fidelity model evaluation in this example. Close examination of~\eqref{eqn:B_thetabar} shows that the controller length scale parameter $\zeta$ is immaterial if $N=1$. The mean discrepancy does not vary with respect to $\z$ since having only one high-fidelity evaluation is not sufficient data to fit the variation. 

To facilitate efficient computation, especially inversion of the shifted linear systems $\alpha \L + \lambda_i \I$, we compute the truncated generalized singular value decomposition of $\vec{D}^{-1}$ in the $\vec{M}$ and $\vec{E}^{-1}$ inner product, which after algebraic manipulations gives the eigenvalue decomposition of $\L^{-1}$. The subsequent linear solves involving $\L$ may be efficiently approximated using the truncated eigenvalue decomposition and the error in the approximation is controlled by the rank, which we take to be $1000$, ensuring a relative truncation error of $\mathcal O(10^{-3})$.

\subsubsection*{Posterior model discrepancy}
\begin{figure}[h]
\centering
\includegraphics[width=0.32\textwidth]{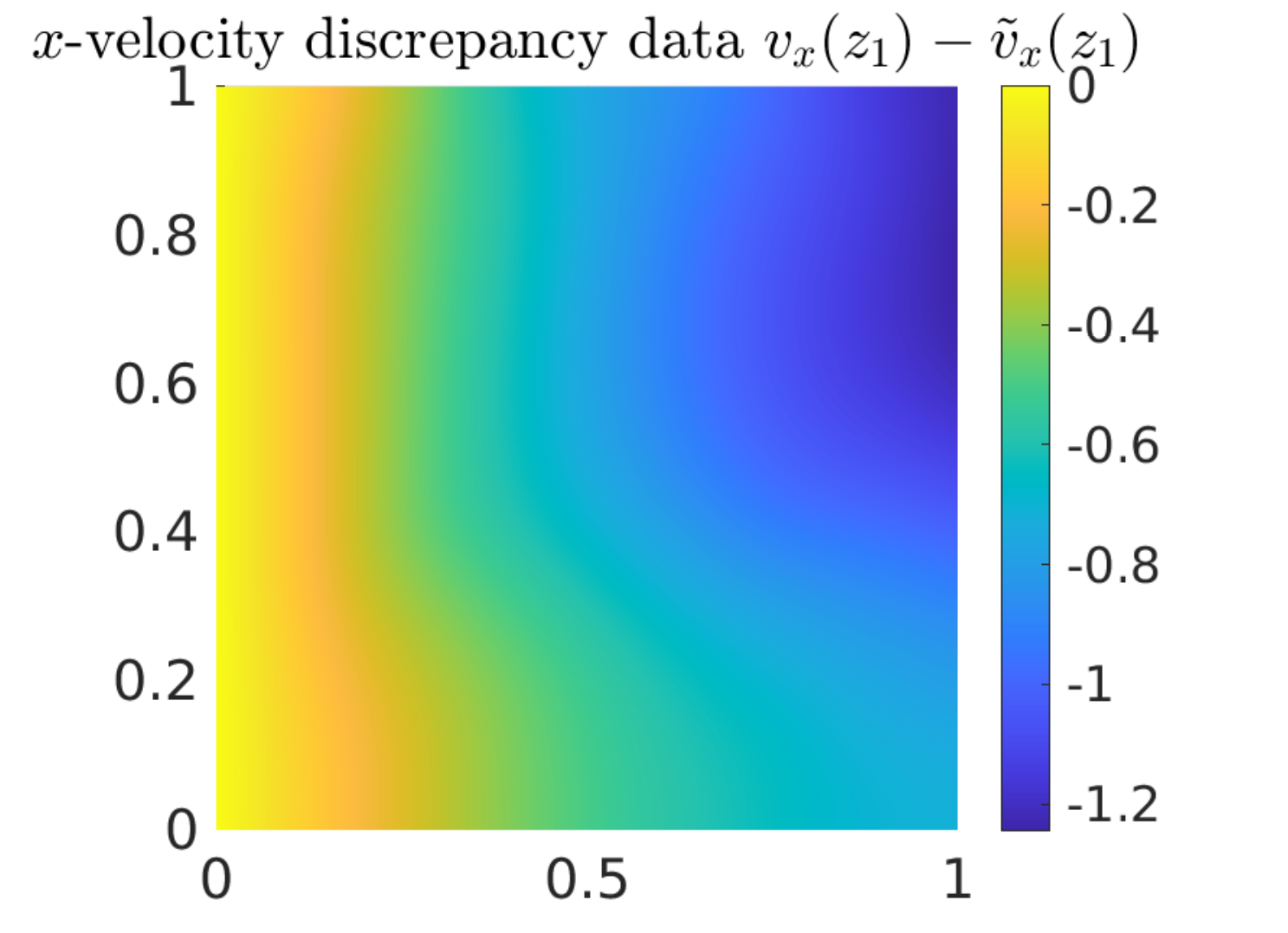}
\includegraphics[width=0.32\textwidth]{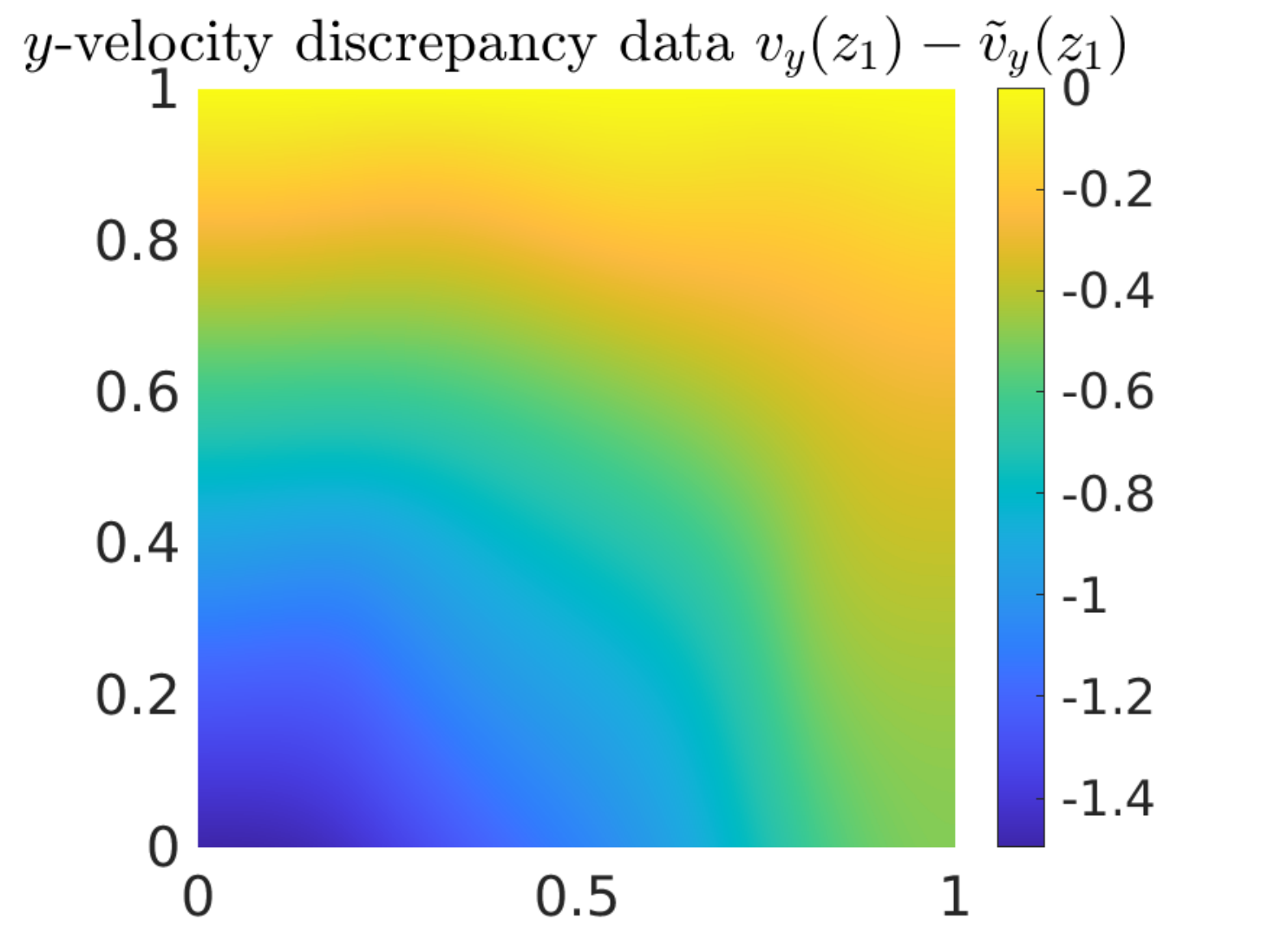}
 \includegraphics[width=0.32\textwidth]{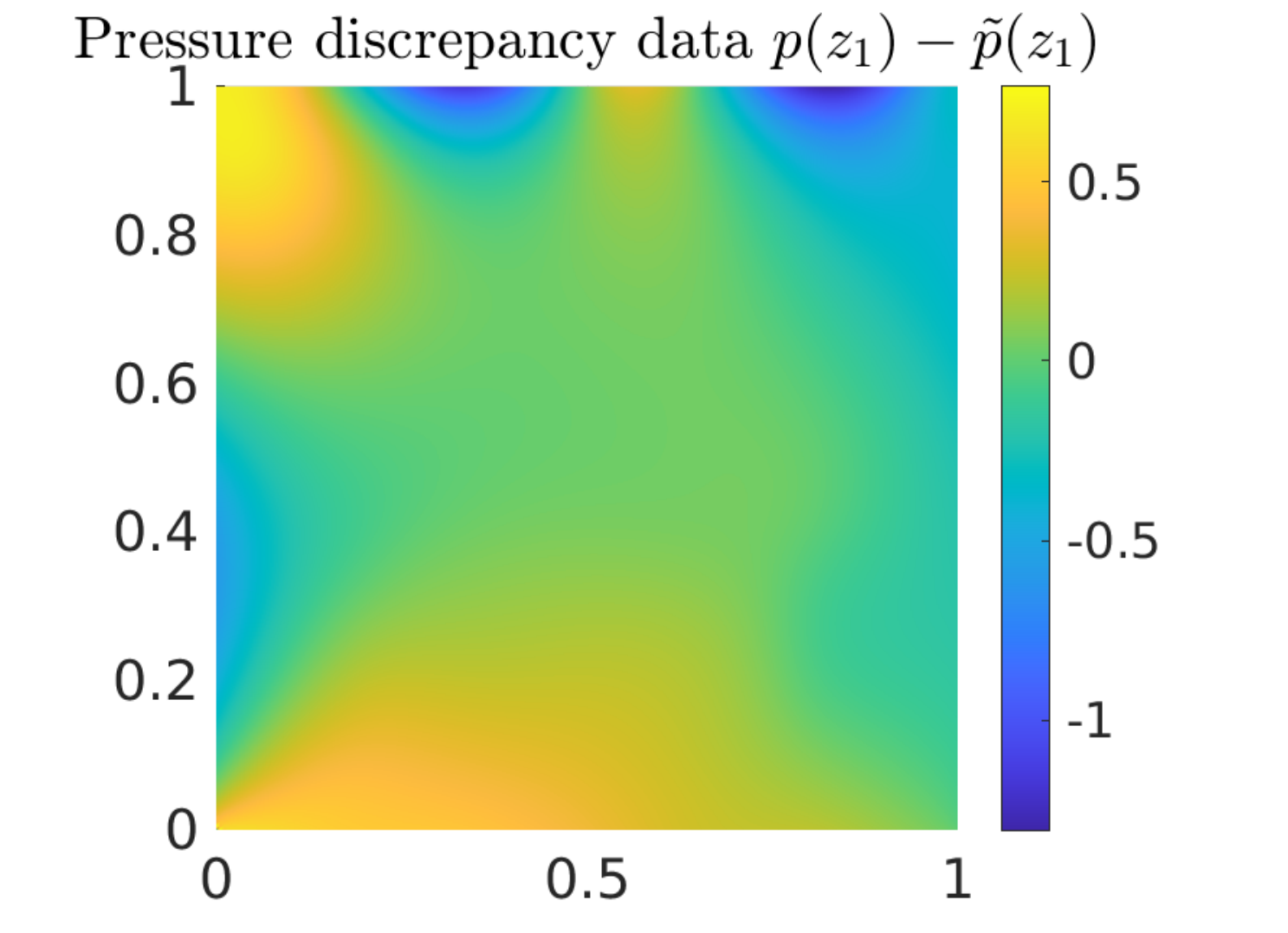} \\
\includegraphics[width=0.32\textwidth]{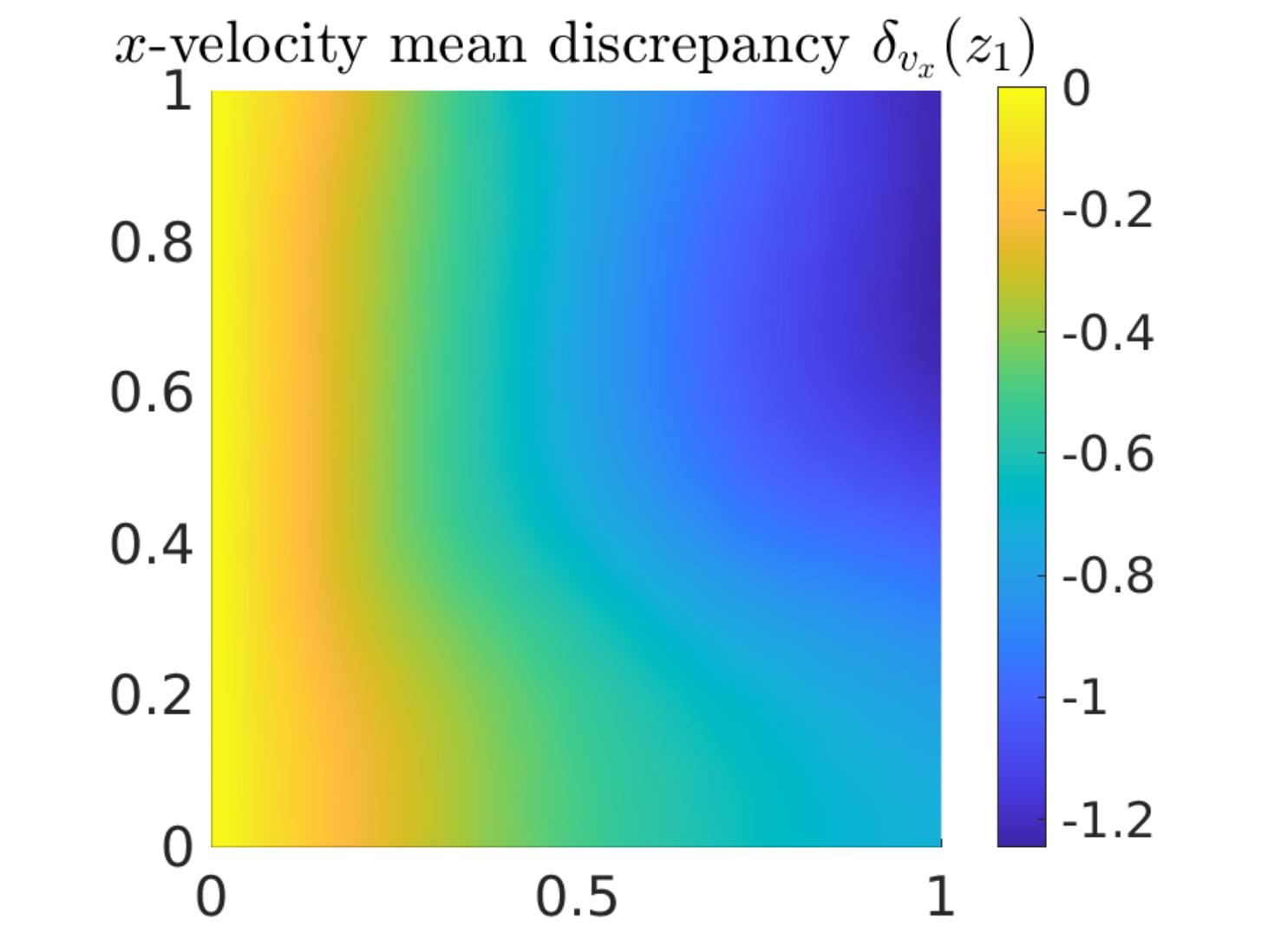}
\includegraphics[width=0.32\textwidth]{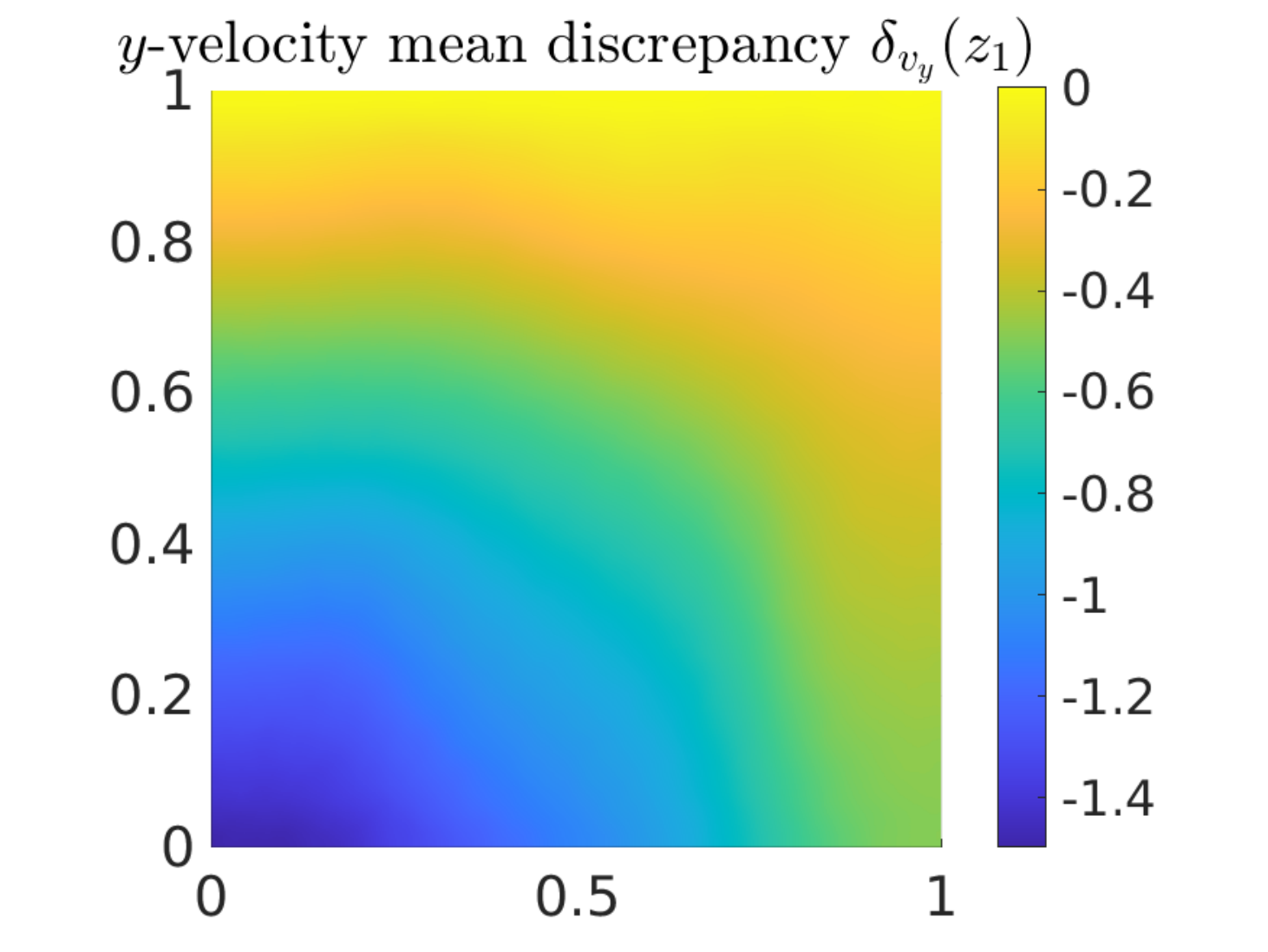}
 \includegraphics[width=0.32\textwidth]{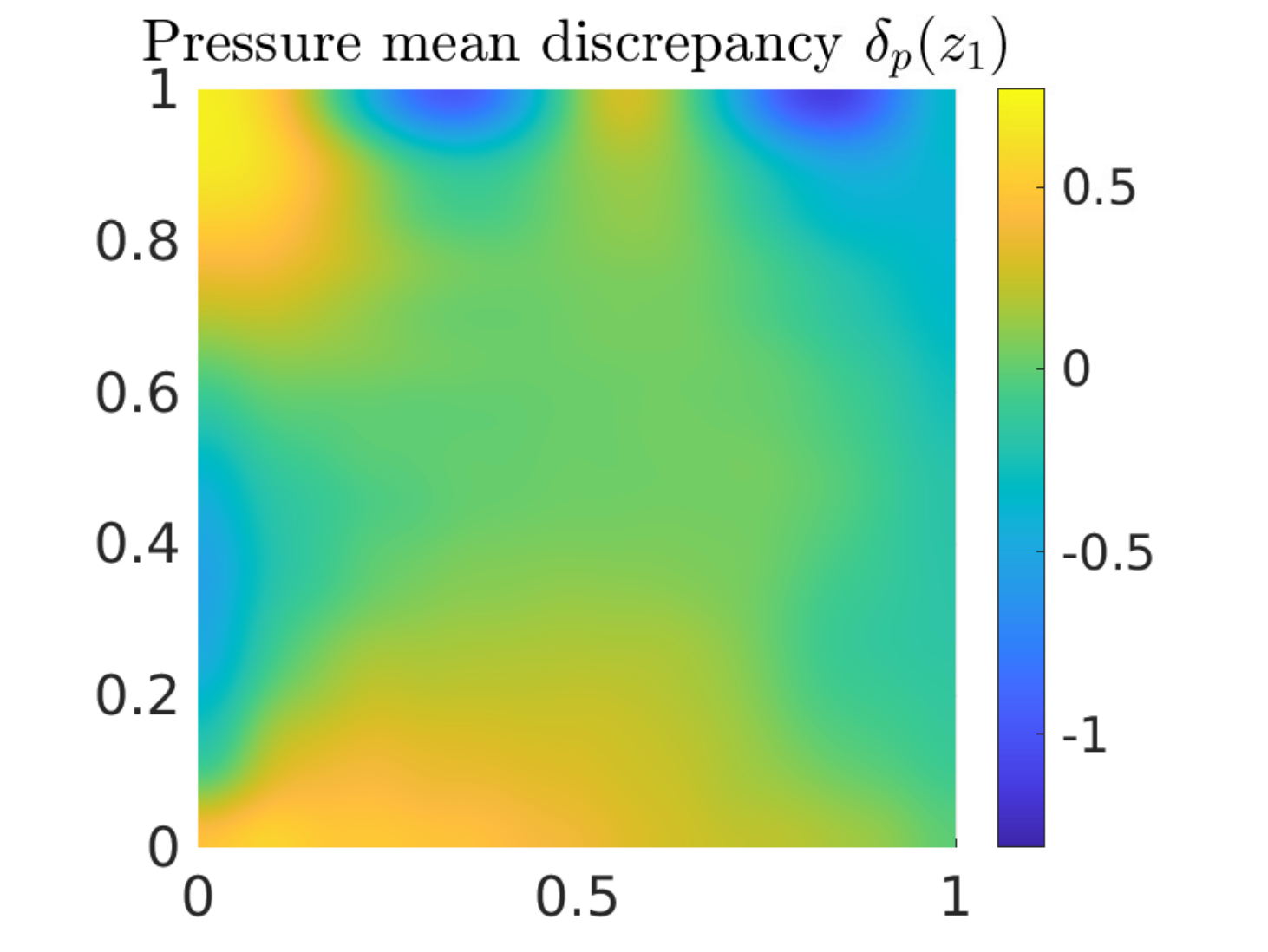}
  \caption{Observed data (top) and mean fit (bottom) for the state discrepancy. The $x$-velocity $v_x$, $y$-velocity $v_y$, and pressure $p$ are plotted from left to right.}
  \label{fig:stokes_example_discrepany_fit}
\end{figure}

We evaluate the Naiver-Stokes equation at the optimal Stokes controller ($N=1$) to generate the data pair $(\z_1,\y_1)$ to fit the model discrepancy. With a setting of $\alpha=0.01$, we solve the Bayesian inverse problem to calibrate $\d(\z,\t)$. Figure~\ref{fig:stokes_example_discrepany_fit} displays the posterior mean of $\d$ evaluated at the controller $\z_1$ and shows a good data fit.

\subsubsection*{Optimal solution update}
\begin{figure}[h]
\centering
\includegraphics[width=0.32\textwidth]{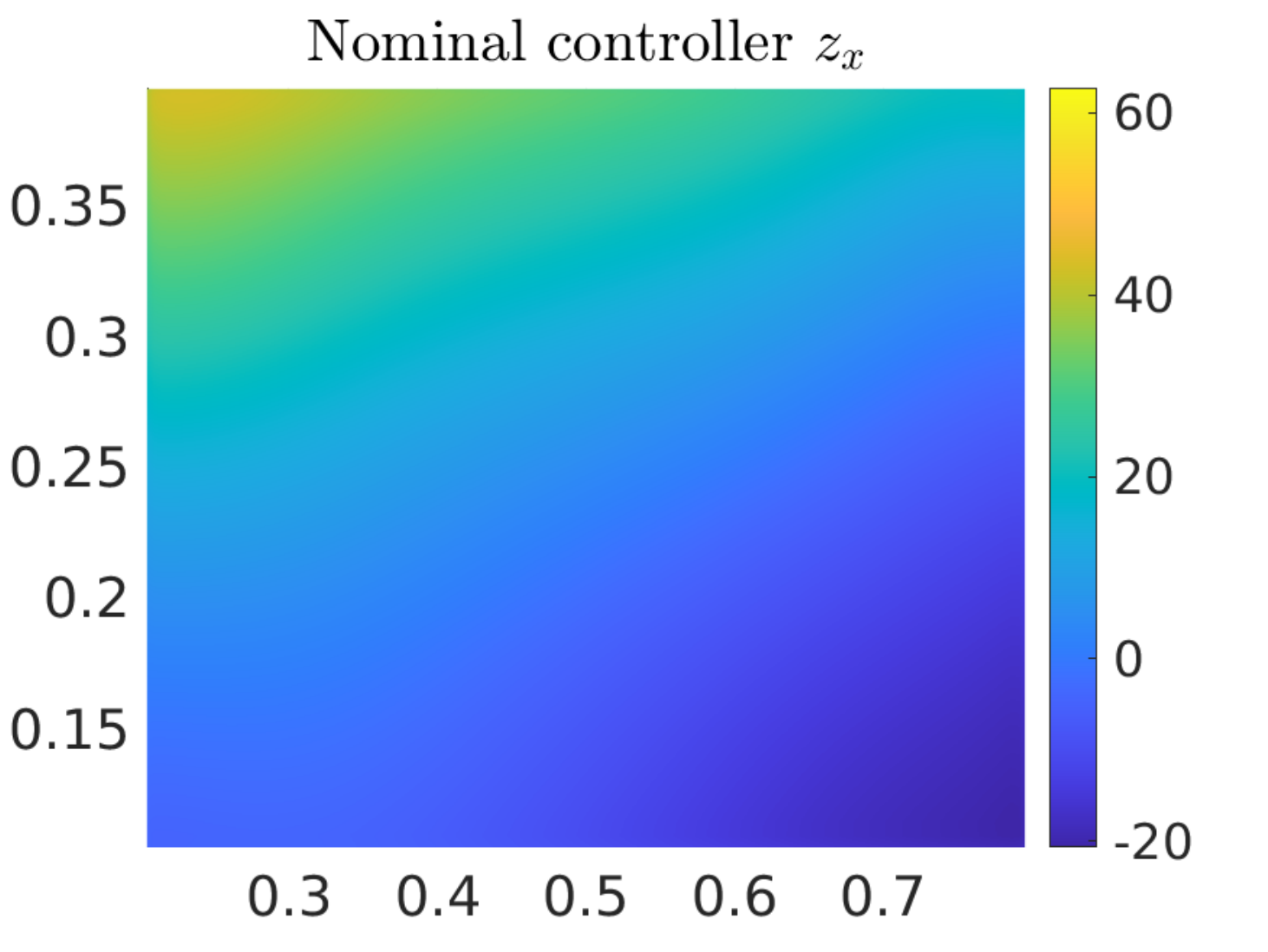}
\includegraphics[width=0.33\textwidth]{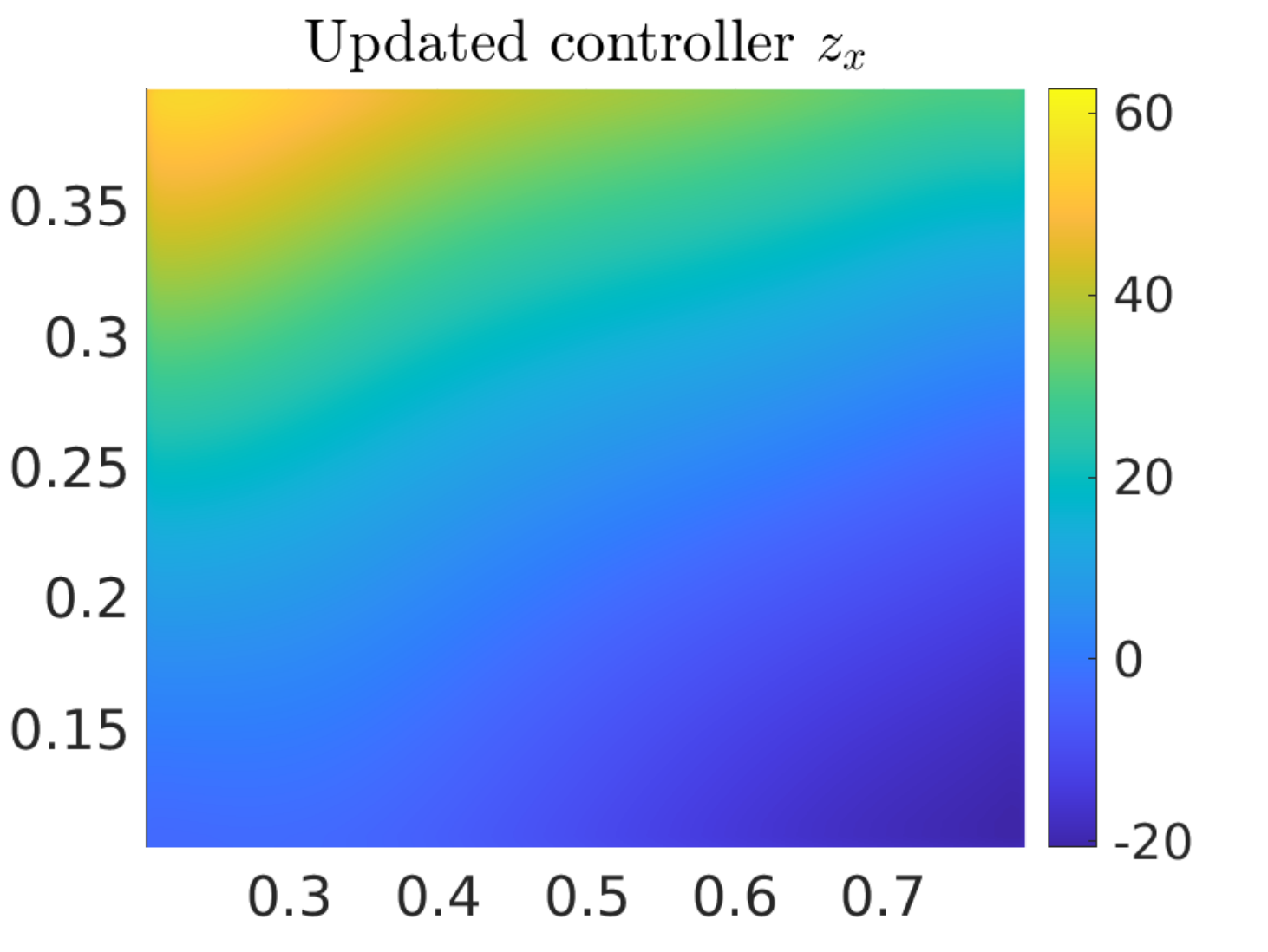}
 \includegraphics[width=0.33\textwidth]{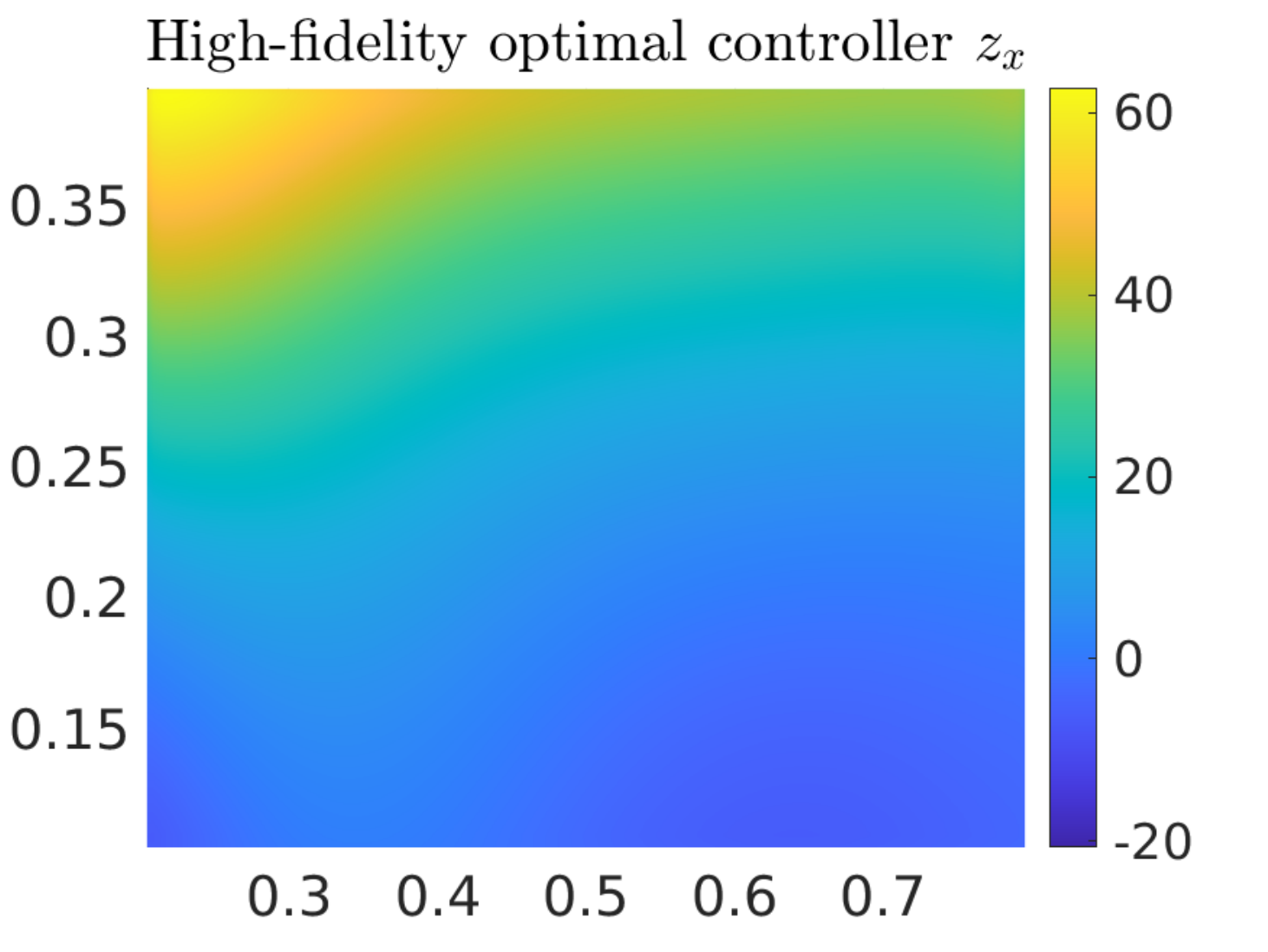} \\
 \hspace{-3mm}
\includegraphics[width=0.33\textwidth]{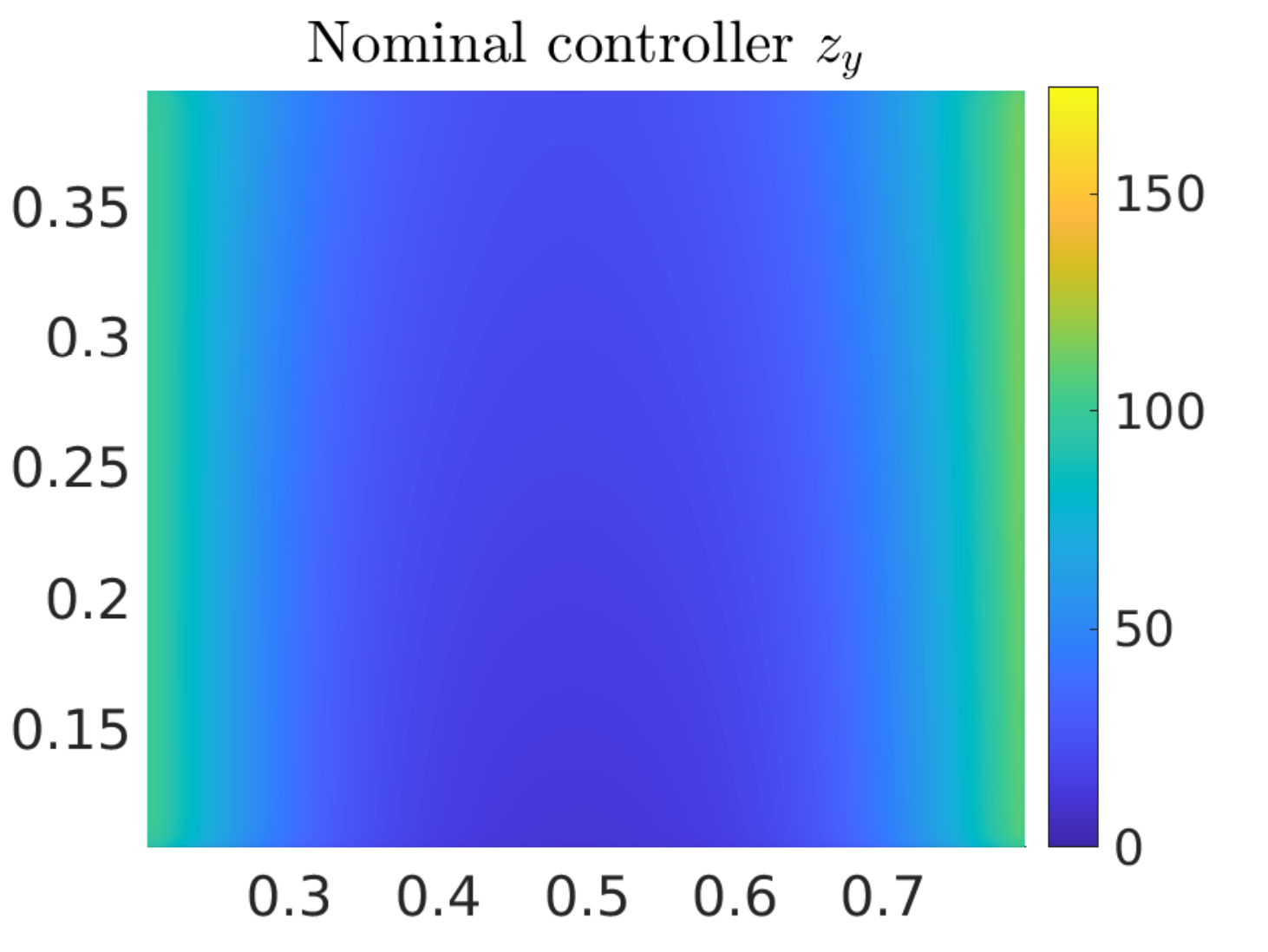}
\includegraphics[width=0.32\textwidth]{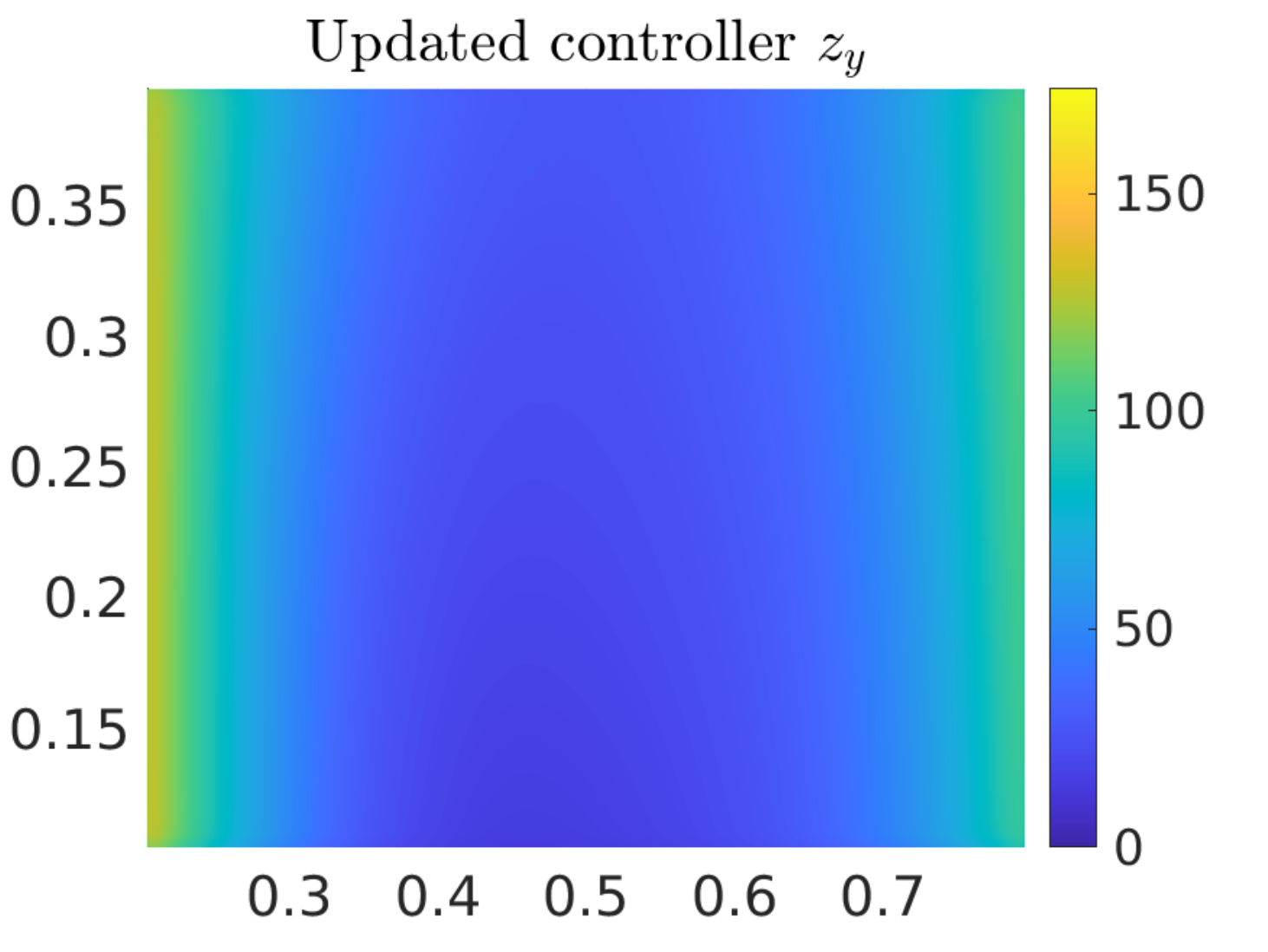}
 \includegraphics[width=0.32\textwidth]{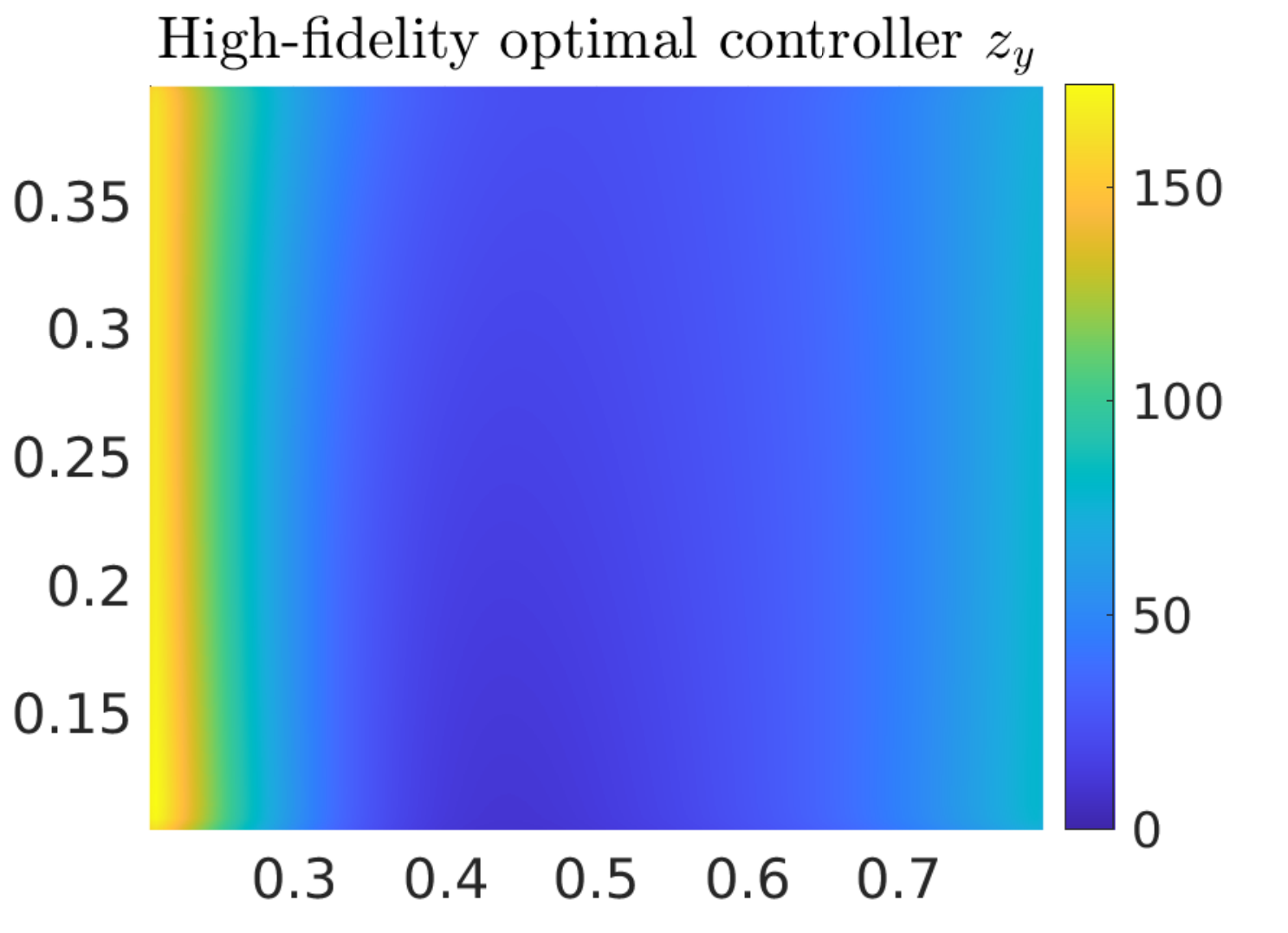} 
  \caption{Controllers acting on the $x$-velocity component (top) and the $y$-velocity component (bottom). The left column is the low-fidelity control solution, the center column is the mean of the posterior optimal solution, and the right column is the high-fidelity optimal solution computed by solving the control problem constrained by the Naiver-Stokes model.}
  \label{fig:stokes_example_controllers}
\end{figure}

Lastly, the mean from the posterior discrepancy is propagated through the post-optimality sensitivity operator to produce the optimal solution update. To illustrate the benefit of the controller update, Figure~\ref{fig:stokes_example_controllers} displays the optimal controller generated by solving the Stokes problem, the updated optimal controller, and the optimal controller generated by solving the control problem constrained by the Naiver-Stokes equation. Although our assumption is that the latter controller is generally not available in practice, we compute and display it here for comparison. Figure~\ref{fig:stokes_example_controllers} demonstrates that calibrating the discrepancy with $N=1$ high-fidelity forward solve significantly improves the controller. 

\begin{figure}[h]
\centering
\includegraphics[width=0.32\textwidth]{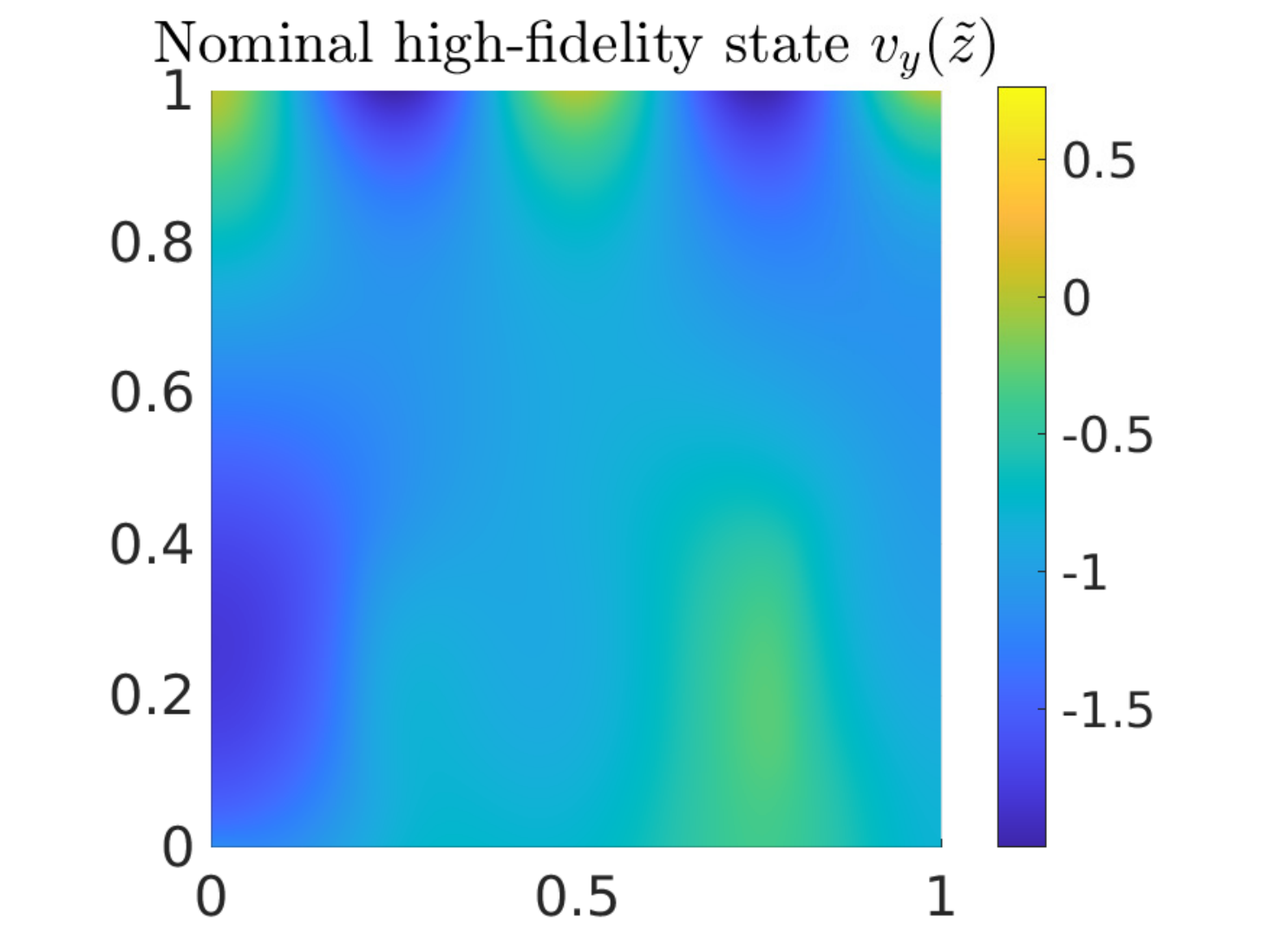}
\includegraphics[width=0.32\textwidth]{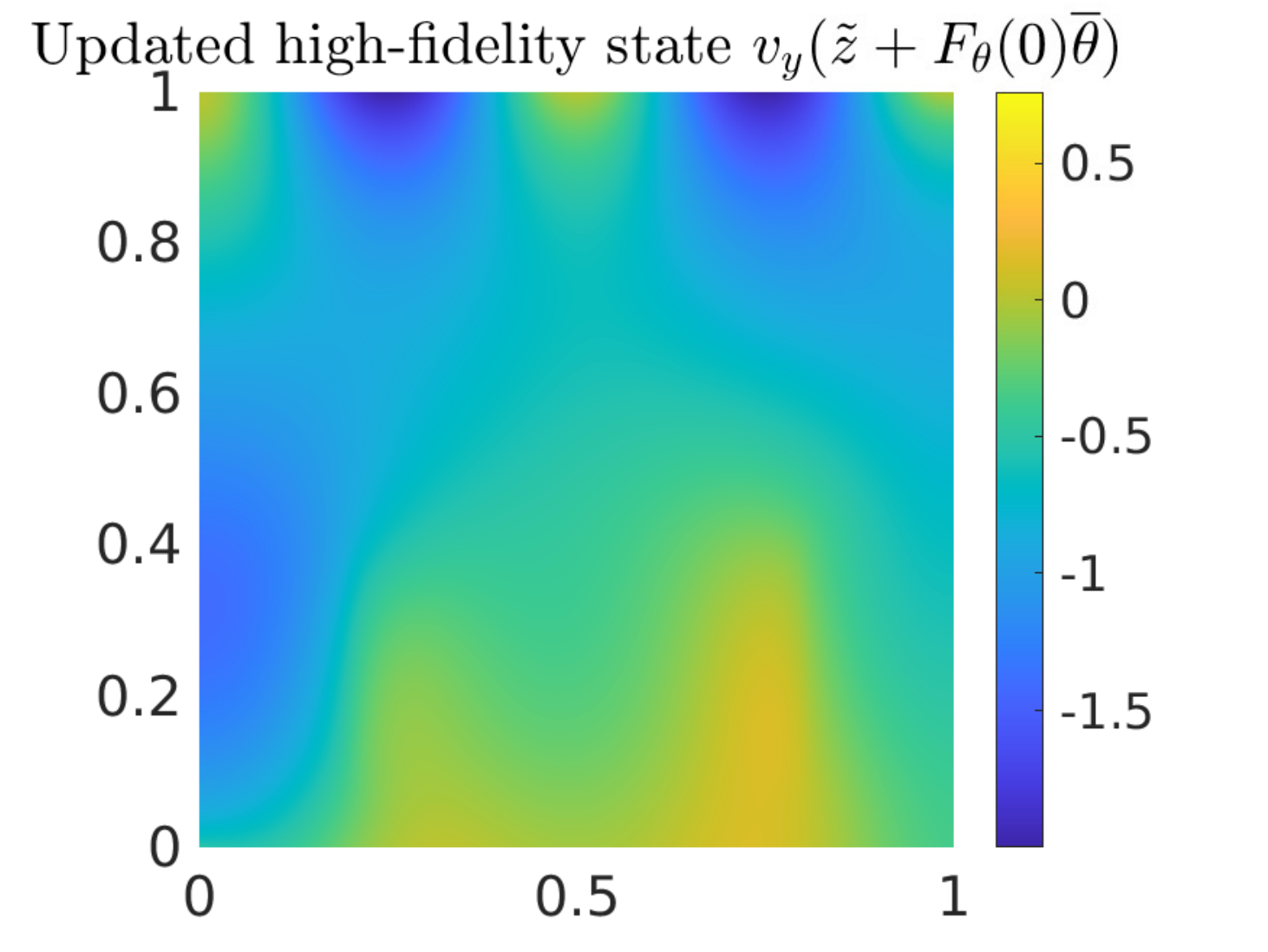}
\includegraphics[width=0.32\textwidth]{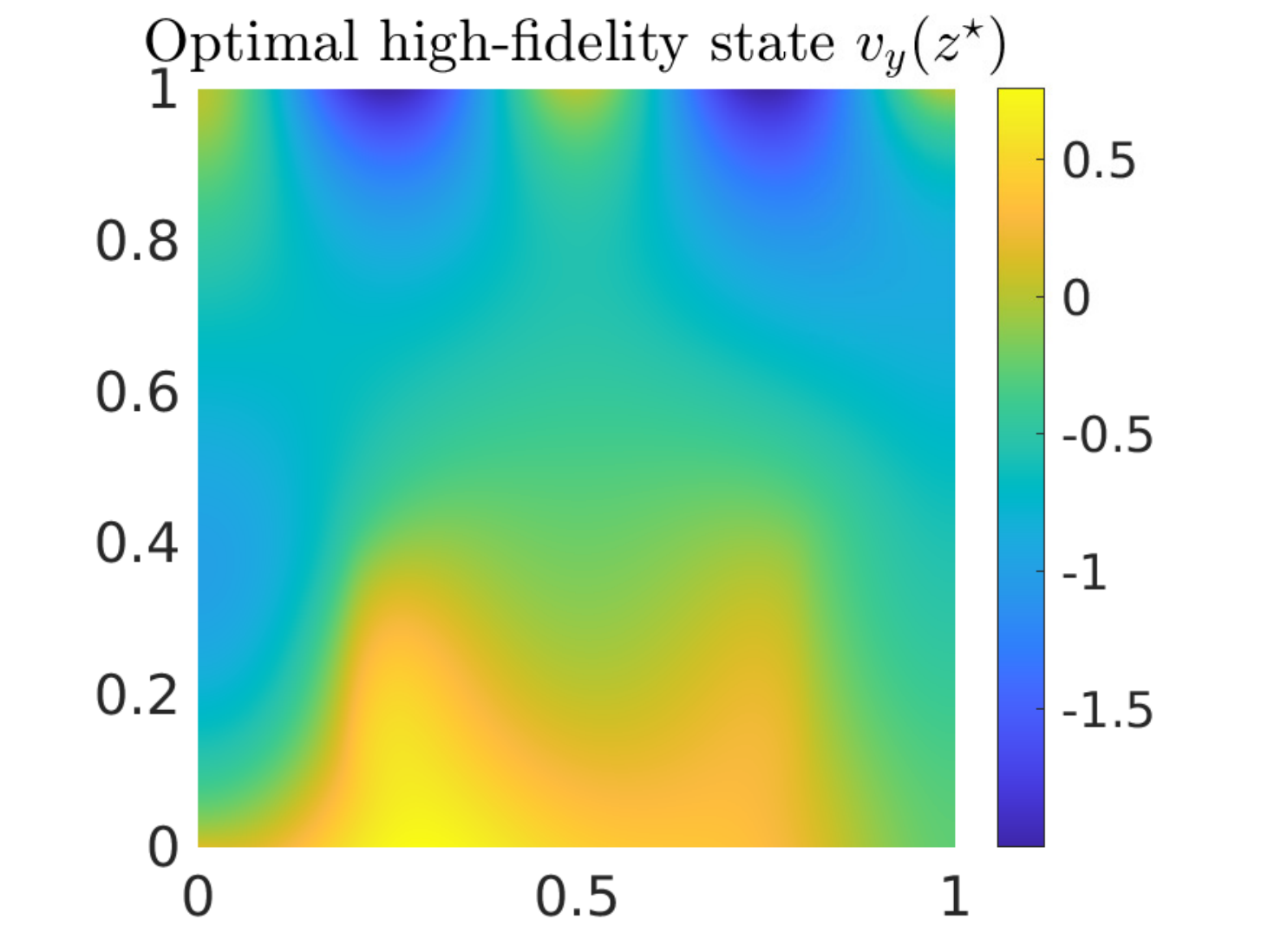}
  \caption{Vertical velocity state solutions corresponding to: Naiver-Stokes solve with low-fidelity control solution (left), Naiver-Stokes solve with updated control solution (center), and Naiver-Stokes with high-fidelity control solution (right).}
  \label{fig:stokes_example_updated_state}
\end{figure}

To further quantify the benefit of the controller update, Figure~\ref{fig:stokes_example_updated_state} displays the vertical velocity $v_y$ of the Naiver-Stokes solution when evaluated at the Stokes optimal controller (left), the updated optimal controller (center), and the Naiver-Stokes optimal controller (right). The performance improvement is significant given that it was achieved using only $N=1$ Naiver-Stokes solve to update the Stokes solution. The extent of the improvement is quantified by the value of the objective function, evaluated using the Naiver-Stokes model, for each of the controllers (Table~\ref{tab:stokes_objective_values}).

\begin{table}[!ht]
\centering
\begin{tabular}{|c|c|}
\hline 
Controller  & Objective Function Value \\
\hline 
$\tilde{\z}$ & $.265$ \\
\hline 
$\tilde{\z} + \nabla_{\t} \F(\vec{0})\overline{\t}$ & $.113$\\
\hline 
$\z^\star$ & $.077$ \\
\hline 
\end{tabular}
\caption{Value of the objective function $J(S(z),z)$, where $S(z)$ is the Naiver-Stokes solution operator.}
\label{tab:stokes_objective_values}
\end{table}

Having the viscosity $\mu=0.5$ in this example ensures that the Reynold number is small. This justifies the Stokes approximation, though Figure~\ref{fig:stokes_example_discrepany_fit} demonstrates that the difference in the Stokes and Naiver-Stokes state solutions is non-trivial. We explored larger Reynolds numbers to stress test the approach. The optimal solution update was less successful in these cases. Nonetheless, we emphasize the significant optimal solution improvement in a nontrivial $\mu=0.5$ case. Future work will explore approaches to quantify how accurate the low-fidelity model must be, and improving on the optimal solution updating strategy when it is not accurate enough. 

\section{Conclusion} \label{sec:conclusion}
It is common that high-fidelity models are too computationally intensive and/or the code structures are not conducive to intrusive implementations of efficient optimization algorithms, for instance, it may be difficult to implement adjoints to enable efficient derivative computations. In such cases, current practice is to construct low-fidelity models through either simplifications of the equations from physics assumptions, or reduced order modeling techniques to approximate the high-fidelity model at a lower computational cost. The low-fidelity model is used to constrain an optimization problem and its solution is assessed by evaluating the high-fidelity model a small number of times. However, these high-fidelity solves may be used for more than assessing the quality of the optimal solution. This article provides a framework that systematically uses the high-fidelity data to improve the optimal solution, without requiring any additional access to the high-fidelity solver. 

Our proposed approach is computationally scalable as a result of judiciously manipulating linear algebra to enable closed form expressions for the posterior. These efficiencies are possible by taking a wholistic perspective that considers every aspect of the problem formulation (discrepancy representation, Gaussian prior and likelihood, post-optimality sensitivity operator, and the resulting matrix factorizations) to develop a mathematically rigorous and computationally advantageous approach. Although the fundamental assumption is a linear approximation of the model discrepancy to optimal solution mapping, we demonstrate considerable improvements for nonlinear PDE-constrained optimization problems with high-dimensional optimization variables.

Since we constrain our analysis in a neighborhood of the low-fidelity solution, there is no guarantee that the proposed approach will converge to the solution of the high-fidelity optimization problem. If we have limited high-fidelity model evaluations, as in our numerical results, then it is unrealistic to expect convergence. This is a sacrifice made to afford an approach that is pragmatic for large-scale applications where high-fidelity simulation data may take hours or days to acquire and high-fidelity derivatives are not accessible. Future work will explore an iterative approach that uses high-fidelity data to update the solution and then takes a new linear approximation about the updated solution. This may extend the applicability of the approach for problems where high-fidelity data is limited but may come in batches of evaluations. 

There are many potential applications of our analysis including, but not limited too, (1) using models in three spatial dimensions to improve optimal solutions computed using models in two spatial dimensions, (2) using multi-scale models to improve solutions computed using homogenized models, (3) using coupled systems to improve optimal solutions computed using only a subset of the system, and (4) using high-fidelity data from controlled experiments in place of a high-fidelity model to improve design or control strategies. Other mathematical questions regarding how and when to query the high-fidelity model also remain and are the topic of ongoing research.

\section*{Acknowledgements}
This paper describes objective technical
results and analysis. Any subjective views or opinions that might be
expressed in the paper do not necessarily represent the views of the
U.S. Department of Energy or the United States Government. Sandia
National Laboratories is a multimission laboratory managed and
operated by National Technology and Engineering Solutions of Sandia
LLC, a wholly owned subsidiary of Honeywell International, Inc., for
the U.S. Department of Energy's National Nuclear Security
Administration under contract DE-NA-0003525. SAND2022-14144 O.

\bibliographystyle{elsarticle-harv} 
\bibliography{dasco}

\end{document}